\documentclass[preprint,5p,times,twocolumn]{elsarticle}
\usepackage{amsmath,amssymb,amsfonts,amsthm}
\usepackage{algorithm}
\usepackage{algorithmic}
\usepackage{graphicx}
\usepackage{textcomp}
\usepackage{dsfont}
\usepackage{framed}
\usepackage{multicol} % Multiple columns environment
\usepackage{nomencl} % Nomenclature package
\makenomenclature
\renewcommand*\nompreamble{\begin{multicols}{2}}
\renewcommand*\nompostamble{\end{multicols}}

\usepackage{array}
\usepackage[caption=false,font=normalsize,labelfont=sf,textfont=sf]{subfig}
\usepackage{textcomp}
\usepackage{stfloats}
\usepackage{url}
\usepackage{verbatim}
\def\BibTeX{{\rm B\kern-.05em{\sc i\kern-.025em b}\kern-.08em
		T\kern-.1667em\lower.7ex\hbox{E}\kern-.125emX}}
\usepackage{balance}
\usepackage{marginnote}
\usepackage{siunitx}
\usepackage{booktabs}
\usepackage[dvipsnames]{xcolor}
\usepackage{multirow}
\usepackage[nolist,nohyperlinks]{acronym}
\usepackage{pbalance}
\usepackage{flushend}
\usepackage[subtle]{savetrees}
\usepackage{natbib}

\newif\ifmargincomments %A quick way of turning off margin comments for, say, arXiv submission
\margincommentstrue

\ifmargincomments
\newcommand{\rwmargin}[2]{#1}%\newcommand{\rwmargin}[2]{{\color{blue}#1}\marginpar{\color{blue}\footnotesize #2}}

\else
\newcommand{\rwmargin}[2]{#1}

\fi

\newtheorem{theorem}{Theorem}[section]

\newcommand{\flexbrac}[1]{\if\relax\detokenize{#1}\relax \else (#1) \fi}
\newcommand{\flexcomma}[1]{\if\relax\detokenize{#1}\relax \else ,#1 \fi}

\acrodef{TCO}[TCO]{Total Cost of Ownership}
\acrodef{BEVs}[BEVs]{{Battery Electric Vehicles}}
\acrodef{EM}[EM]{{Electric Motor}}
\acrodef{EoS}[EoS]{{Economy-of-Scale}}

\begin{document}

\begin{frontmatter}

 \title{Concurrent Design Optimization of Powertrain Component Modules \\ in a Family of Electric Vehicles~\tnoteref{label1}}
 \tnotetext[label1]{This publication is part of the project NEON with project number 17628 of the research program Crossover which is (partly) financed by the Dutch Research Council (NWO).}
 \author{Maurizio~Clemente\corref{correspondingAuthor}\fnref{addressTUe}}
 \ead{m.clemente@tue.nl}
 \cortext[correspondingAuthor]{Corresponding author.}
	\author{Mauro~Salazar\fnref{addressTUe}}
	\author{Theo~Hofman\fnref{addressTUe}}
	\affiliation[addressTUe]{organization={Department of Mechanical Engineering,
		Eindhoven University of Technology},%Department and Organization
	addressline={Groene Loper}, 
	city={Eindhoven},
	postcode={5600MB}, 
	country={The Netherlands}}
%		\thanks{This publication is part of the project NEON with project number 17628 of the research program Crossover which is (partly) financed by the Dutch Research Council (NWO).}% <-this % stops a space
%	\thanks{$^{1}$, 
%		{\tt\small \{m.clemente,m.r.u.salazar,t.hofman\}@tue.nl}}%and hence cost of operation
	
	\begin{abstract}
		We present a modeling and optimization framework to design powertrains for a family of electric vehicles, focusing on the concurrent sizing of their motors and batteries.
		Whilst tailoring these component modules to each individual vehicle type can minimize energy consumption, it can result in high production costs due to the variety of component modules to be realized for the family of vehicles, driving the Total Costs of Ownership (TCO) high.
		Against this backdrop, we explore modularity and standardization strategies whereby we jointly design unique motor and battery modules to be installed in all the vehicles in the family, using a different number of these modules when needed.
		Such an approach results in higher production volumes of the same component module, entailing significantly lower manufacturing costs due to Economy-of-Scale (EoS) effects, and hence a potentially lower TCO for the family of vehicles.
		To solve the resulting ``one-size-fits-all'' problem, we instantiate a nested framework consisting of an inner convex optimization routine which jointly optimizes the modules' sizes and the powertrain operation of the entire family, for given driving cycles and modules' multiplicities.
		Likewise, we devise an outer loop comparing each configuration to identify the minimum-TCO solution with global optimality guarantees.
		Finally, we showcase our framework on a case study for the Tesla vehicle family in a benchmark design problem, considering the Model~S, Model~3, Model~X, and Model~Y.
		Our results show that, compared to an individually tailored design, the application of our concurrent design optimization framework achieves a significant reduction of the production costs for a minimal increase in operational costs, ultimately lowering the family TCO in the benchmark design problem by 3.5\%.
		Moreover, our concurrent design optimization methodology can reduce the TCO by up to 17\% for the market conditions considered in our sensitivity study.
	\end{abstract}

%	This approach significantly reduces production costs by leveraging economy-of-scale effects (higher volumes entail spreading overhead costs on larger numbers of items) at the expense of higher overall energy consumption, as the components are now designed to equip multiple vehicles in a one-size-fits-all fashion.
	
	
\end{frontmatter}

\newtheoremstyle{break}% name
{}%         Space above, empty = `usual value'
{}%         Space below
{\itshape}% Body font
{}%         Indent amount (empty = no indent, \parindent = para indent)
{\bfseries}% Thm head font
{.}%        Punctuation after thm head
{\newline}% Space after thm head: \newline = linebreak
{}%         Thm head spec
\theoremstyle{break}

\newtheorem{problem}{\textbf{Problem}}

% Input sections
\section{Introduction}\label{sec:introduction}

%	we aim to enable the access of a larger share of the population to electric vehicles (EVs) while still accounting for their ever-increasing range of needs.
%	 , reckoning the diversification opportunities for vehicle manufacturers and components suppliers.
%	For this reason, we 
%	capturing the impact of modules' changing size and multiplicity on the mechanical power demand, performance, and energy consumption of the vehicles 
%	we identify a convex model of the powertrain,
%(SOCP)
%	First,.
%	Second, we develop a cost model to estimate the TCO of vehicles, taking into account the electricity price, production volumes, vehicle type, multiplicity and size of the components in the production and operation costs.
%	Third, every possible combination of the number of modules in the family is framed as a second-order conic program , jointly optimizing the shared modules' size.
%	Finally, we explore all the combinations 
%		 Yet such an advantage could come at the cost of a lower vehicle efficiency compared to individually tailored BEVs.
%		In an effort to speed up the transition to sustainable mobility, we aim to reduce the EVs' total cost of ownership (TCO) leveraging product-family design strategies, like modularity and standardization, to trigger economy-of-scale effects and drive down the production costs of components.

The wider diffusion of \ac{BEVs} is considered a key objective in the transition to more sustainable mobility~\cite{IEA2021}.
As a matter of fact, \ac{BEVs} are important tools in the fight against pollution in cities, lowering the environmental impact by significantly decreasing particulate matter~\cite{McTurk2022} and CO$_2$ emissions~\cite{IEA2020}. %allies\tools
Despite the efforts in promoting policies oriented to the purchase and use of \ac{BEVs}~\cite{EGR21}, they still account for a small share of the total amount of vehicles due to a relatively higher price when compared to their conventional, fossil-fuel-powered counterparts~\cite{IEA2022}.
The considerable upfront cost, combined with supply shortages affecting manufacturing~\cite{Paoli2022}, prevents the adoption of this technology by a larger share of the general public~\cite{IEA2020,Energy2021}.
%The reduction in manufacturing cost translates into a lower cost of the vehicle, which solves one of the major issues preventing a larger uptake of electric mobility~\cite{IEA2020}. 
In order to accelerate the transition to electric mobility, the \ac{TCO} of \ac{BEVs} must be reduced.
Many companies have adopted design strategies aimed at developing remarkably efficient tailor-made designs~\cite{Lightyear2022}. %minimize operational costs by 
However, these design choices imply considerable production costs, due to the highly specific components that are individually tailored to each vehicle model.

\begin{figure}[t]
	\centering
	\includegraphics[width=\columnwidth]{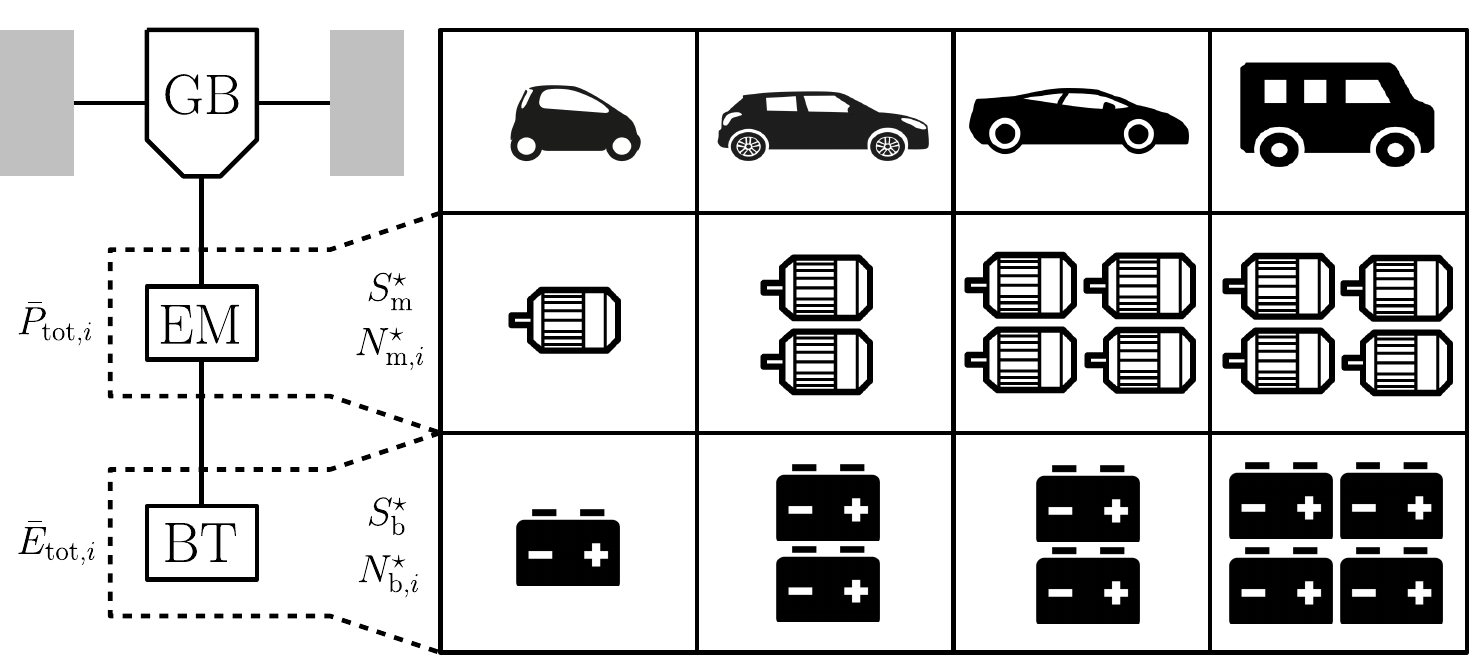}
	\caption{Family of battery electric vehicles designed leveraging a shared modular powertrain. Every vehicle type is equipped with $N_{\mathrm{m},i}^\star$ motor modules (EM) of size $S_{\mathrm{m}}^\star$ and $N_{\mathrm{b},i}^\star$ battery modules of size $S_{\mathrm{b}}^\star$, for a total maximum motor power of $\overline{P}_{\mathrm{tot},i}$ and maximum energy capacity of $\overline{E}_{\mathrm{tot},i}$.}
	\label{fig:modular}
\end{figure}

\begin{figure}[t]
	\centering
	\includegraphics[width=0.88\columnwidth]{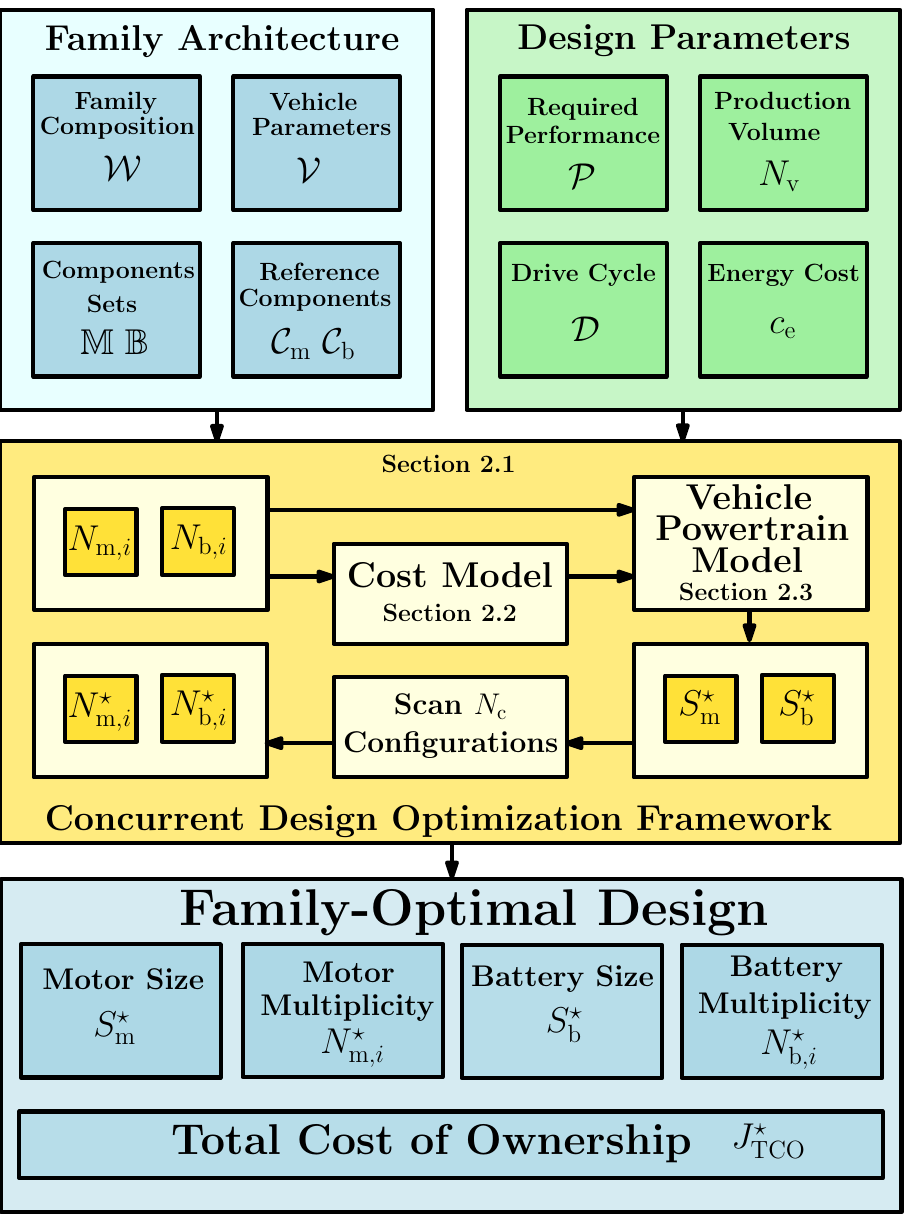}
	\caption{Concurrent design optimization methodology diagram.}
	\label{fig:concurrent}
\end{figure}

Conversely, it is possible to reduce production costs by leveraging \ac{EoS} strategies, whereas increasing the amount of identical components produced reduces their specific production cost.
However, this kind of approach requires large production volumes of the same item, involving a critical manufacturing ramp-up~\cite{MediniPierneEtAl2020} and hindering diversification, pivotal to satisfying the wide range of customer needs.
We aim to address these issues by investigating a module-based product family design approach (Fig.~\ref{fig:modular}) to take full advantage of the \ac{EoS} strategies, substantially lowering the vehicles' \ac{TCO} without limiting product diversification. %Many companies have adopted design strategies aimed at developing remarkably efficient tailor-made designs
Nevertheless, the mere identification of a feasible design of a standardized powertrain module for multiple vehicles is a challenging task, which may result in excessive performance deterioration in the vehicles~\cite{AMO2014}.
For this reason, we resort to the application of numerical optimization to find the shared size and number of modules that every vehicle is equipped with (\textit{multiplicity}), jointly optimized to minimize the \ac{TCO} of the family instead of being individually tailored for every vehicle.
%whilst accounting for the changing module’s size as well as multiplicity. 
Against this backdrop, this paper presents the concurrent design optimization framework shown in Fig.~\ref{fig:concurrent} to jointly optimize the size and multiplicity of shared battery and \ac{EM} modules with the operations of a family of \ac{BEVs}.%, minimizing the overall \ac{TCO}.

% (\textit{variants})
\emph{Related Literature:}
This paper pertains to two main research lines: multi-product design and powertrain design optimization.
Multi-product design is a design strategy leveraging commonality, modularity and standardization among products to reduce components' costs, provide operational and logistical advantages in part sourcing, and quality control~\cite{JiaoSimpsonEtAl2007}.
Moreover, it fosters the development and upgrade of differentiated products efficiently, increases flexibility and responsiveness in manufacturing processes~\cite{RobertsonUlrich1998}, and generates substantial savings in research, testing, interface design, and integration~\cite{OttoHoelttae-OttoEtAl2016}.
Multi-product design consists of two main processes: platform design and product family design. The former entails the development of a product platform (consisting of  a common product architecture, shared physical components and processes) from the company objectives and market research. The latter accounts for the design of products starting from the common platform.
Industrial players have widely studied and employed these methodologies for different products to provide cost-effective variety and customization~\cite{SimpsonSiddiqueEtAl2006}: from the Sony Walkman~\cite{SandersonUzumeri1995,SandersonUzumeri1997} to aircraft~\cite{Sabbagh1996,RothwellGardiner1990} and spacecraft~\cite{CaffreySimpsonEtAl2022}.
%``In the past decade, there has been a flurry of research activity in the engineering design community to develop methods and tools to facilitate product platforms and product family design ".
The automotive sector was perhaps one of the fields where multi-product design achieved the best results.
Industrial giants like Volkswagen~\cite{Bremmer1999,Bremmer2000} saved billions of dollars per year producing some of the most successful automotive platforms. % due to the ``greater flexibility between plants and increased plant usage"~\cite{Muffatto1999}.
In fact, ``in the 1990s, automotive manufacturers that employed a platform-based product development approach gained a 5.1\% market share per year while those that did not, lost 2.2\%" \cite{NobeokaCusumano1997}.
%Bremmer supports this claim by considering that ``Volkswagen saved an estimated \$1.5 billion per year in development and capital costs using platforms, and they produced three of the six automotive platforms that successfully achieved production volumes over one million in 1999".
However, the application of these strategies to electric vehicles has been delayed due to the smaller production volumes compared to conventional vehicles~\cite{IEA2021}, the lack of standardization~\cite{PereirinhaTrovao}, and the challenges in design due to the delayed ripening technology.

The second research stream concerns powertrain design optimization.
This discipline has been thoroughly studied over the years, producing a considerable number of methods and applications~\cite{GuzzellaSciarretta2007}.
In the last decade, powertrain design in \ac{BEVs} has witnessed significant developments due to the blossoming of new system-level design optimization strategies~\cite{SilvasHofmanEtAl2016,WangZhouEtAl2022} involving higher complexity and larger design spaces~\cite{DongZhaoEtAl2023,KwonLimEtAl2023}.
There are many examples of design strategies focusing on the joint optimization of the system and controls, aiming at reducing energy consumption~\cite{HofmanSalazar2020,RadrizzaniRivaEtAl2023,VerbruggenRangarajanEtAl2019,BorsboomFahdzyanaEtAl2021,CarlosDaSilvaKefsiEtAl2023,BorsboomSalazarEtAl2022,ZhangQinEtAl2020}, thus driving down the operation costs of the vehicle.
Other authors analyzed the trade-off between the components' cost and the vehicle's energy efficiency~\cite{RibauSilvaEtAl2014}, concluding that they are conflicting objectives in the optimization.
Nevertheless, none of these methods considers the trade-off between the energy efficiency of a vehicle-tailored design against the significant manufacturing cost reduction prompted by designing shared modules serving a whole family of vehicles.
%Nonetheless, these methods only account for the design optimization of a single vehicle, not considering the trade-off between the energy efficiency of a vehicle-tailored design against the drastic cost reduction prompted by the \ac{EoS} strategy of designing shared modules serving a whole family of vehicles.
The application of product family design concepts to hybrid electric vehicles can be found in one recent paper~\cite{Anselma2022}, where the author proposed a tool to compare different topologies and sizes from a predefined set to minimize the costs for the manufacturer.
 %``Consistent advantages may be achieved from the Original Equipment Manufacturer (OEM) perspective in terms of streamlining the manufacturing process and reducing the number of different parts managed, thus laying the foundations for easing development and mass production processes of more complex powertrains such as the electrified ones''
However, the algorithm can only compare a limited number of size values from predefined sets for each component in the powertrain, and the cost model does not consider the savings owing to the \ac{EoS} effects.
%Finally, it only focuses on HEV, not considering Battery Electric Vehicle (BEV) powertrains.
%to the best of the authors’ knowledge, 
%As a matter of fact, the combined application of module-based product family concepts and vehicles optimization has not been studied extensively.
%Fellini and Kokkolaras~\cite{FelliniKokkolarasEtAl2002,FelliniKokkolarasEtAl2002b} used optimization for making commonality decisions while controlling individual performance in a family of cars and developed a sensitivity-based commonality strategy for family products of mild variation. Yet their application concerns only automotive body structures.
%Therefore, we identified a research gap in the application of product-family methodologies to the optimal design of powertrain modules for a family of BEV.% concerning the cost model e algoritm?
Although the concept of modular design, the use of an optimization approach, and convex optimization techniques are not new in themselves, to the best of the authors' knowledge, there are no product family design optimization frameworks for electric powertrains capturing the impact of the \ac{EoS} with global optimality guarantees.

%which employs a product family design approach
\emph{Statement of Contribution:}
In this paper, we reconsider the traditional vehicle-tailored optimal design methodology for BEVs in favor of a uniquely sized family-optimal module design, seeking the best compromise between the cost reduction induced by the \ac{EoS} and energy efficiency.
In particular, we instantiate a nested framework jointly optimizing the battery and \ac{EM} modules' size and multiplicity with vehicles' operations to minimize the family \ac{TCO}, explicitly accounting for the effects of production volumes.
%In furtherance of this goal, we develop:
To this end, we develop:
\begin{itemize}
	\item A cost model estimating the vehicle's \ac{TCO} by capturing the influence of production volumes, energy cost, and modules' size and multiplicity;
	\item A vehicle and powertrain operation model, including scalable \ac{EM} and battery modules, taking into account changing modules' sizing and multiplicities;
	\item A low-level convex optimization routine, jointly optimizing the sizing of the modules with the vehicles' operations to minimize the family \ac{TCO};
	\item A high-level algorithm exploring all the vehicles' multiplicity configurations of the family, each with its own optimized modules' size, to identify the one globally minimizing the overall \ac{TCO}.
\end{itemize}
A preliminary version of this paper~\cite{ClementeSalazarEtAl2022} was presented at the IFAC Symposium on Advances in Automotive Control in 2022.
This extended version includes a broader literature review, considers a transmission in the vehicle model, introduces a cost model explicitly accounting for the effect of production volumes, and allows for the optimization of the modules' multiplicity together with their sizing.
Furthermore, we present a real-world case study considering a leading manufacturer vehicle family, and a sensitivity study on the impact of energy costs and vehicles' production volumes on the optimal design solution.
In this benchmark problem, we assumed the same type fraction for each vehicle to demonstrate our BEV design framework independently of the sales reported by the specific company.
%\textbf{AAAAAAAA}
%\rwmargin{Even though we considered a realistic benchmark problem, tackling the co-design of different vehicle types belonging to the Tesla family, the main focus of this paper is to present a BEV optimal family design methodology. It is beyond the scope of our manuscript to analyze sales, market dynamics, and strategic business decisions.}{R5:13}

% (\ref{subsec:probdef}) (\ref{subsec:VECM}) (\ref{subsec:CM})  (\ref{subsec:OPF}) (\ref{subsec:discussion})

\emph{Organization:} The remainder of this paper is organized as follows: Section \ref{sec:methodology} introduces the concurrent optimization design methodology, describes the vehicle's cost and powertrain models, frames the optimization problem formulation, and discusses the assumptions and limitations of our approach.
In Section \ref{sec:results}, we demonstrate the applications of our framework with a benchmark problem where we identify the optimal modules' size and multiplicity for the Tesla vehicle family.
Furthermore, we produce a sensitivity analysis of the methodology advantages under different energy costs and production volumes.
The conclusions are drawn in Section \ref{sec:conclusions}, together with an outlook on future research.

\section{Methodology}\label{sec:methodology}

\begin{table*}[!t]   
	
	\begin{framed}
		
		\nomenclature{$S_\mathrm{b}$}{Battery Scaling Factor}
		\nomenclature{$S_\mathrm{m}$}{Motor Scaling Factor}
		\nomenclature{$N_\mathrm{b}$}{Battery Multiplicity}
		\nomenclature{$N_\mathrm{m}$}{Motor Multiplicity}
		
		\nomenclature{$N$}{Vehicle Types}
		\nomenclature{$N_\mathrm{v}$}{Vehicle Production Volume}
		
		\nomenclature{$P_\mathrm{0}$}{Constant Motor Loss Coefficient}
		\nomenclature{$\beta$}{Linear Motor Loss Coefficient}
		\nomenclature{$\alpha$}{Quadratic Motor Loss Coefficient}
		\nomenclature{$a_k$}{Linear Battery Loss Coefficient}
		\nomenclature{$b_k$}{Constant Battery Loss Coefficient}
		\nomenclature{$\overline{P}_\mathrm{m,0}$}{Peak Reference Motor Power}
		\nomenclature{$\overline{T}_\mathrm{m,0}$}{Maximum Reference Motor Torque}
		\nomenclature{$\omega_\mathrm{r}$}{Rated Motor Speed}
		\nomenclature{$m_\mathrm{m,0}$}{Reference Motor Mass}
		\nomenclature{$m_\mathrm{b,0}$}{Reference Battery Mass}
		\nomenclature{$\overline{E}_\mathrm{b,0}$}{Full Reference Battery Capacity}
		\nomenclature{$\gamma$}{Gear Ratio}
		\nomenclature{$r_\mathrm{FWD}$}{Front Wheel Drive Regen. Braking Fraction}
		\nomenclature{$r_\mathrm{AWD}$}{All Wheel Drive Regen. Braking Fraction}
		\nomenclature{$\xi$}{State-of-charge}	
		\nomenclature{$w_i$}{Optimization Weight of the Vehicle}	
		\nomenclature{$m_\mathrm{g}$}{Glider Mass}			
		\nomenclature{$m_\mathrm{p}$}{Payload Mass}
		\nomenclature{$m_\mathrm{d}$}{Driver Mass}			
		\nomenclature{$\eta_\mathrm{gb}$}{Gearbox Efficiency}			
		\nomenclature{$\eta_\mathrm{inv}$}{Inverter Efficiency}		
		\nomenclature{$r_\mathrm{w}$}{Wheel Radius}			
		\nomenclature{$c_\mathrm{r}$}{Rolling Resistance Coefficient}
		\nomenclature{$c_\mathrm{r}$}{Drag Coefficient}		
		\nomenclature{$A_\mathrm{f}$}{Frontal Area}			
		\nomenclature{$P_\mathrm{aux}$}{Auxiliaries Power}
		\nomenclature{$g$}{Gravitational Constant}	
			
		\nomenclature{$t_\mathrm{a}$}{Acceleration Time}			
		\nomenclature{$v_\mathrm{f}$}{Final speed after Acceleration}			
		\nomenclature{$v_\mathrm{t}$}{Top Speed}			
		\nomenclature{$d_\mathrm{r}$}{Driving Range}				
		\nomenclature{$\theta$}{Road Slope}		
		\nomenclature{$v_\mathrm{m}$}{Uphill-driving speed}
		
		\nomenclature{$J_\mathrm{TCO}$}{Family Total Cost of Ownership}	
		\nomenclature{$j_\mathrm{TCO}$}{Vehicle Total Cost of Ownership}	
		\nomenclature{$C_\mathrm{a}$}{Vehicle Acquisition Price}	
		\nomenclature{$C_\mathrm{op}$}{Vehicle Operation Cost}
		\nomenclature{$C_\mathrm{p}$}{Vehicle Production Cost}	
		\nomenclature{$k_\mathrm{oh}$}{Overhead Cost Factor}	

		\nomenclature{$E_\mathrm{v}$}{Vehicle Lifetime Energy Consumption}	
		\nomenclature{$c_\mathrm{e}$}{Energy Cost}	
		\nomenclature{$F_\mathrm{v}$}{Distance-specific Energy Consumption}	
		\nomenclature{$d_\mathrm{v,lt}$}{Vehicle Lifetime}	
		\nomenclature{$C_\mathrm{g}$}{Glider Cost}	
		\nomenclature{$C_\mathrm{m}$}{Motor Cost}
		\nomenclature{$C_\mathrm{b}$}{Battery Cost}	
		\nomenclature{$c_\mathrm{m}$}{(Peak) Power-specific Motor Cost}
		\nomenclature{$\lambda_\mathrm{m}$}{Economy-of-scale Motor Coefficients}
		\nomenclature{$c_\mathrm{b}$}{Capacity-specific Battery Cost}
		\nomenclature{$\lambda_\mathrm{b}$}{Economy-of-scale Battery Coefficients}
		
		\nomenclature{$P_\mathrm{v}$}{Driving Cycle Required Power}
		\nomenclature{$P_\mathrm{m}$}{Motor Output Power}							
		\nomenclature{$P_\mathrm{ac}$}{Motor Input Power}	
		\nomenclature{$P_\mathrm{b}$}{Battery Output Power}			
		\nomenclature{$P_\mathrm{i}$}{Battery Internal Power}			
		\nomenclature{$P_\mathrm{sc}$}{Battery Short Circuit Power}
		\nomenclature{$E_\mathrm{b}$}{Battery Module Energy}
		\nomenclature{$E_\mathrm{tot}$}{Overall Battery Energy}					
			
%		\nomenclature{$P_\mathrm{v}$}{Driving Cycle Required Power}			
		\printnomenclature
		
	\end{framed}
	
\end{table*}

In this section, we illustrate our framework in detail. Section \ref{subsec:probdef} describes the concurrent design optimization problem, while Sections \ref{subsec:CM} and \ref{subsec:VECM} characterize the cost and powertrain model, respectively.
Finally, we formalize the optimization problem formulation in Section \ref{subsec:OPF} and discuss the assumptions and limitations of our approach in Section \ref{subsec:discussion}.

%This section presents the ``Concurrent Design Optimization" framework.
%In Section~\ref{subsec:sandn} we introduce the design variables of the problem.
%In Section~\ref{subsec:ptconsumption} we describe the vehicle's energy consumption minimization problem while in Section~\ref{subsec:costmodel} we set forth the cost model of the vehicle, including economy-of-scale effects, to estimate the production cost, selling price, and cost of operation of every vehicle in the family.
%Finally, section~\ref{subsec:multiplicity} formalizes the overarching framework determining the optimal modules' configuration in the family by confronting all the possible combinations to find the lowest TCO.

\subsection{Problem Definition: Concurrent Design Optimization}\label{subsec:probdef}

\begin{figure}[t]
	\centering
	\includegraphics[width=\columnwidth]{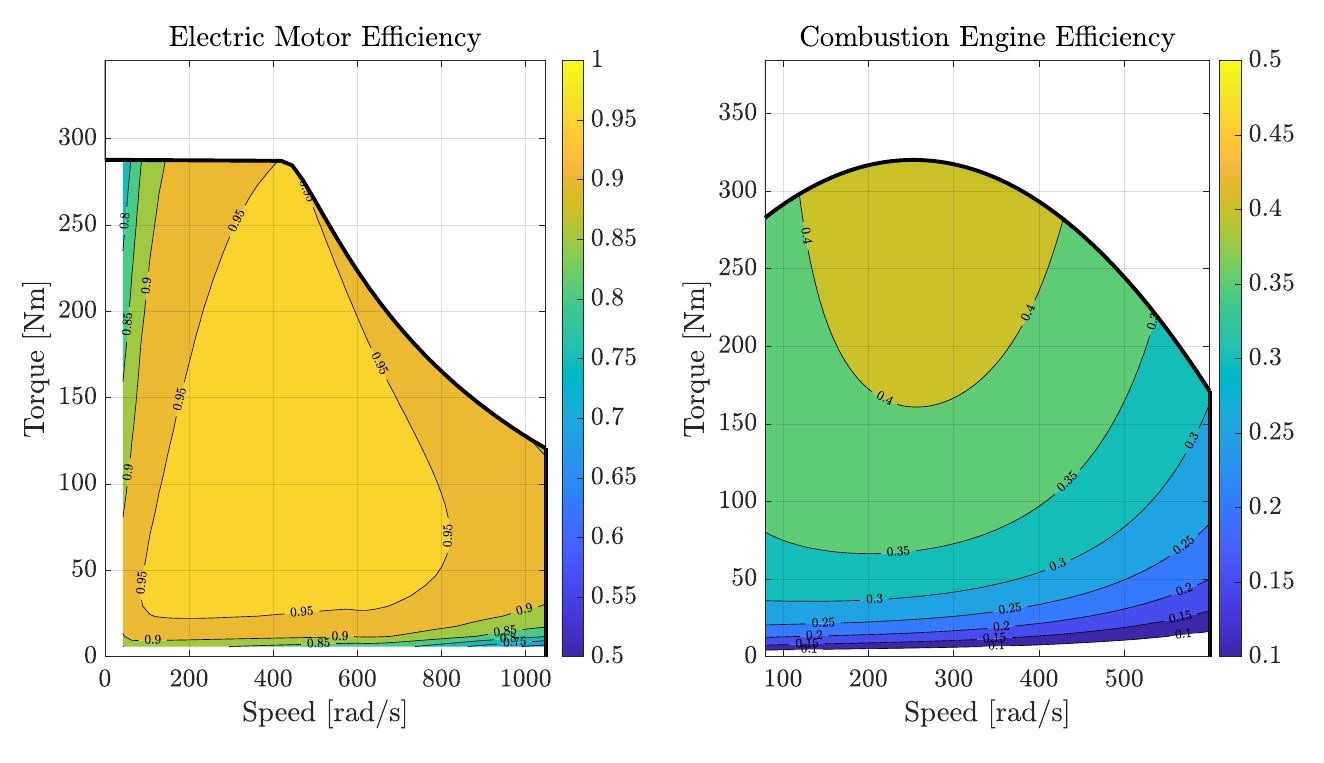}
	\caption{Efficiency map of an electric motor (left) compared with an internal combustion engine (right). Data from \cite{Solipuram2024}. }
	\label{fig:ICEEV}
\end{figure}

The traditional powertrain design approach consists of tailoring the components' sizes to each individual vehicle.
When considering electric vehicles, the large high-efficiency area in an \ac{EM} torque map allows for merely small power losses in many different operating points, as opposed to the conventional internal combustion engine's relatively tight sweet spot region (Fig.~\ref{fig:ICEEV}).
This feature enables differently-sized vehicles to share the same motor model with limited drawbacks, which can be further reduced when employing a modular design, taking advantage of the topologies enabled by the technology.
Furthermore, the battery is an intrinsically modular component and can be easily manufactured in packs of the family-optimal size, satisfying the energy requirements with the optimal number of modules and fostering customization and diversification (range-extended variants).
The great advantage of using the same module type for many different vehicles lies in the \ac{EoS} effects triggered by the increased production volumes.
Compared to a vehicle-tailored design, where every vehicle features components specifically optimized for its design and operations, in concurrent design the manufacturer can amortize overhead costs on a larger number of items, reducing their specific cost.
Although this approach enables battery and motor swapping business, our design methodology does not focus on this feature, prioritizing the cost reduction induced by standardization and modularization.
%Thus, more vehicles sharing the same components produce more significant absolute savings in manufacturing costs.
However, due to the asymptotic behavior of the phenomenon, usually described as the ``law of diminishing returns", an increase in production volumes will bear relatively smaller advantages in module cost every time.
Despite the module cost reduction, the standardization of the components could lead to underwhelming performance or excessively high energy consumption, negatively impacting the end-user of the vehicle.
%For these reasons, it's important to consider the \ac{TCO} in the analysis, including the cost of operations together with the price of acquisition of the vehicle, which, considering fixed percentual overhead costs, accounts for the production costs of the manufacturer.
For this reason, we exploit numerical optimization techniques to strike the optimal balance between energy efficiency and \ac{EoS} cost reduction, taking into account all the constraints from each of the vehicles in the family simultaneously, to avoid performance degradation.
We devise a bi-level, nested algorithm exploring every configuration of the modules in the family and leveraging convex properties to rapidly converge to the globally optimal solution.
Subsequently, we compare the optimally sized configurations to identify the lowest \ac{TCO} while retaining information on the other sub-optimal solutions.
%and minimize the \ac{TCO} of the family, accounting for both the end-users and manufacturer through its two constituents: the cost of operations and the vehicle acquisition price, which can be related to the manufacturing cost~\cite{KoenigNicolettiEtAl2021}.
%the cost of operations and the selling price of the vehicle, which can be linked to the manufacturing cost, considering both the end-user and the manufacturer
% including the cost of operations together with the price of acquisition of the vehicle, which, considering fixed percentual overhead costs~\cite{KoenigNicolettiEtAl2021}, accounts for the production costs of the manufacturer.
%This way, we can leverage this methodology to design a family of \ac{BEVs} for the.

Hereby, we present the design variables of the problem: the motor scaling factor $S_\mathrm{m}$, the battery scaling factor $S_\mathrm{b}$, and the respective number of modules (\textit{multiplicity}) $N_{\mathrm{m},i}$ and $N_{\mathrm{b},i}$.
The subscript $i$ indicates that the quantity differs from one vehicle type to the other, as opposed to their common sizing.
We denote the minimum and maximum value of a variable $X$ by introducing a line below and above its symbol, $\underline{X}$ and $\overline{X}$, respectively.
We devise our model using a reference motor and a reference battery for the identification of the parameters.
In order to preserve a clean notation, these parameters are included in the sets of reference components $\mathcal{C_\mathrm{m}}$ and $\mathcal{C_\mathrm{b}}$
\begin{equation*}\label{eq:Cm}
	\mathcal{C}_\mathrm{m} \coloneqq \left\{ P_\mathrm{0}, \beta, \alpha, m_{\mathrm{m,o}}, \overline{P}_{\mathrm{m,o}}, \omega_{\mathrm{r}}, \gamma, r_{\mathrm{FWD}}, r_{\mathrm{AWD}}  \right\},
\end{equation*}
\begin{equation*}\label{eq:Cb}
	\mathcal{C}_\mathrm{b} \coloneqq \left\{ a_k, b_k, m_{\mathrm{b,o}}, \overline{E}_\mathrm{b,o}, \underline{\xi}, \overline{\xi} \right\},
\end{equation*}
where the motor loss coefficients $ P_\mathrm{0} $, $ \beta $, and $ \alpha $ are dependent on the motor speed $\omega$ and subject to identification. Similarly, the battery loss coefficients $a_k$ and $b_k$, $k\in\{1,...,K\}$, are identified from a reference battery and used to determine the short-circuit power $P_{\mathrm{sc}}$ of the battery, a measure of efficiency described in Section \ref{subsubsec:BT}. %and appendix \ref{app:OCP}.
The constants $m_{\mathrm{m,o}}$ and $m_{\mathrm{b,o}}$ are the reference masses of the motor and battery, respectively.
The parameters $\underline{\xi}$ and $\overline{\xi}$ are the minimum and maximum state-of-charge operational limits, while $r_{\mathrm{FWD}}$ and $r_{\mathrm{AWD}}$ are the regenerative braking fraction for a Front-Wheel Drive and an All-Wheel Drive.
Finally, $\gamma$ is the gear ratio and $\omega_{\mathrm{r}}$ is the speed at which the maximum torque and maximum power curves intersect, also called rated speed.
Consistently, we assume that the vehicle's maximum output power of the motor(s) $\overline{P}_{\mathrm{tot},i}$ and the maximum energy of the battery pack $\overline{E}_{\mathrm{tot},i}$ are obtained by linearly scaling the reference motor's maximum output power $\overline{P}_{\mathrm{m,o}}$ and  reference battery maximum energy capacity $\overline{E}_\mathrm{b,o}$, taking into account the modules' multiplicity
\begin{equation}
	\overline{P}_{\mathrm{tot},i} = S_\mathrm{m} \cdot N_{\mathrm{m},i} \cdot \overline{P}_{\mathrm{m,o}},
	\label{eq:Sm}
\end{equation}
\begin{equation}
	\overline{E}_{\mathrm{tot},i} = S_\mathrm{b} \cdot N_{\mathrm{b},i} \cdot \overline{E}_\mathrm{b,o}.
	\label{eq:Sb}
\end{equation}
Nonetheless, these approximations are only valid in the range of scales
\begin{equation}\label{eq:Smcommon}
	S_\mathrm{m} \in \left[\underline{S}_\mathrm{m}, \overline{S}_\mathrm{m}\right] \subseteq \mathbb{R_+},
\end{equation}
\begin{equation}\label{eq:Sbcommon}
	S_\mathrm{b} \in \left[\underline{S}_\mathrm{b}, \overline{S}_\mathrm{b}\right] \subseteq \mathbb{R_+},
\end{equation}
while the possible multiplicities of the modules in the vehicles are listed in the sets $\mathbb{M},\mathbb{B}$
\begin{equation}\label{eq:Nm}
	N_{\mathrm{m},i} \in \mathbb{M} \coloneqq \left\{\underline{N}_\mathrm{m}, ..., \overline{N}_\mathrm{m}\right\} \subseteq \mathbb{N_+},
\end{equation}
\begin{equation}\label{eq:Nb}
	N_{\mathrm{b},i} \in \mathbb{B} \coloneqq \left\{\underline{N}_\mathrm{b}, ..., \overline{N}_\mathrm{b}\right\} \subseteq \mathbb{N_+}.
\end{equation}
Moreover, we introduce the Family Composition $\mathcal{W}$ and Vehicle Parameters $\mathcal{V}$ sets from the Family Architecture block of Fig. \ref{fig:concurrent},
\begin{equation*}\label{eq:F}
	\mathcal{W} \coloneqq \left\{ w_1, w_2, ... , w_i \right\},
\end{equation*}
\begin{equation*}\label{eq:F}
	\mathcal{V}_i \coloneqq \left\{ m_{\mathrm{g},i}, {m}_{\mathrm{p},i}, {m}_{\mathrm{d},i}, \eta_\mathrm{gb}, \eta_\mathrm{inv}, r_{\mathrm{w},i}, c_{\mathrm{r},i}, c_{\mathrm{d},i}, A_{\mathrm{f},i}, P_\mathrm{aux}  \right\} \subseteq \mathcal{V},
\end{equation*}
where the $w_i$ are the fractions of the total number of vehicles $N_{\mathrm{v}}$ of the $i$-th type, while the set $\mathcal{V}$ includes $N$ subsets $\mathcal{V}_i$, one for each vehicle type considered.
Each subset $\mathcal{V}_i$ contains, in respective order, the glider mass $m_{\mathrm{g},i}$, the payload mass ${m}_{\mathrm{p},i}$, the gearbox efficiency $\eta_\mathrm{gb}$, the inverter efficiency $\eta_\mathrm{inv}$, the wheel radius $r_{\mathrm{w},i}$, the frontal area $A_{\mathrm{f},i}$, and the rolling resistance $c_{\mathrm{r},i}$ and aerodynamic drag $c_{\mathrm{d},i}$ coefficients.
Furthermore, we include the required performance that the vehicle has to satisfy in the set $\mathcal{P}$, with its $N$ subsets
\begin{equation*}\label{eq:P}
	\mathcal{P}_i \coloneqq \left\{ \overline{t}_{\mathrm{a},i}, v_{\mathrm{f}}, \underline{v}_{\mathrm{t},i}, \underline{d}_{\mathrm{r},i}, \underline{v}_{\mathrm{m},i}, \underline{\theta}_i \right\} \in \mathcal{P},
\end{equation*}
where $\overline{t}_{\mathrm{a},i}$ is the maximum time needed to accelerate from $0$ to $v_{\mathrm{f}}$, $\underline{v}_{\mathrm{t},i}$ the minimum required top speed, $\underline{d}_{\mathrm{r},i}$ the minimum required range, $\underline{v}_{\mathrm{m},i}$ is the minimum speed at which the vehicle shall be able to drive facing a slope of at least $\underline{\theta}_i$.

Ultimately, using our concurrent design optimization framework, we are able to identify the family-optimal sizes $ S_\mathrm{m}^\star, S_\mathrm{b}^\star$ and multiplicities $N_{\mathrm{m},i}^\star, N_{\mathrm{b},i}^\star $  of the modules minimizing the family \ac{TCO} $J_{\mathrm{TCO}}^\star$ by striking the optimal compromise between production and operational costs.

\subsection{Cost Model}\label{subsec:CM}

The majority of authors who developed a cost model for \ac{BEVs} predominantly analyze the decline in costs over time.
These analyses are based on technology-development predictions, which are affected by a high degree of uncertainty~\cite{HoekstraVijayashankarEtAl2017,SimeuKim2018,SimeuKimEtAl2021,FriesKerlerEtAl2017,KochhanFuchsEtAl2014,BrookerGonderEtAl2015}.
Conversely, we introduce a more robust model based on the relation between cost of motors and batteries and production volumes instead of uncertain technology development forecasts. For this reason, we identified the asymptotic dependence of the components’ cost on the production volume from manufacturing data~\cite{Anderman2019,KoenigNicolettiEtAl2021,Lipman1999,Cuenca1995,PhilippotAlvarezEtAl2019,BEAN21,SimeuKim2018, SimeuKimEtAl2021,CoxBauerEtAl2020}, and validated with data from the market~\cite{Tesla2023,I-MIEV2023,EVDatabase2022,Twingo2023,DACIA2023,OPEL2023,Peugeot2023,KiaNiro2023,Polestar2023}.
%We also show the numerical values of the identified coefficients in Table \ref{tab:costpar} in the Result Section, among other parameters of the cost model.
We introduce an asymptotic dependence of the components' cost on the production volume, based on manufacturing data, to account for the impact of \ac{EoS} effects.
%The benefits of \ac{EoS} strategies are confidential information, often available merely to experts in the field, who can effectively sell their know-how via technical reports or consulting services.}{@Mauro ok here or to be moved somehow to literature review section?}
%With this cost model, we aim to describe the impact of components' production volumes on the manufacturing cost, providing an estimation tool for researchers and industry, enabling further analyses on the most profitable conditions to apply \ac{EoS} strategies.
We display the data used in Table \ref{tab:data} of \ref{app:CM}: First, we portray the asymptotic curves resulting from Eq. \eqref{cm} and \eqref{cb} from Section \ref{subsubsec:mc} (Motor Cost) and Section \ref{subsubsec:bc} (Battery Cost) in Figures \ref{fig:Masy} (motor) and \ref{fig:Basy} (battery pack).
Then, we validate the results of the cost model in Table \ref{tab:costcomparison}, providing a visual representation in Fig. \ref{fig:costcomparison}.
However, we do not claim any correlation of production cost with sales, as this relation depends on the analysis of complicated customer-market dynamics which are beyond the scope of this paper.
In our concurrent design optimization methodology, the description of these effects enables analyzing the trade-off in \ac{TCO} between a tailored design, where the components have larger production costs and higher efficiencies, and a family-shared design, with smaller costs and lower efficiencies.
For the purpose of this paper, we assume fixed or vehicle-class dependent costs such as insurance, maintenance, and taxes to be independent of the decision variables in the optimization framework~\cite{KoenigNicolettiEtAl2021}.
While these additional costs could be readily included in the model, the minor differences in motor and battery sizing compared to the (individual) vehicle-tailored design do not significantly affect the relative outcome of our study.
For this reason, in a first approximation, we disregard them in our analysis.
Hereby, we define the family \ac{TCO} as the sum of the \ac{TCO}s of all the vehicles in the family $j_{\mathrm{TCO},i}$:
\begin{equation}\label{eq:TCO}
	J_{\mathrm{TCO}} = \sum_{i=1}^{N} w_i \cdot j_{\mathrm{TCO},i}.
\end{equation}
Each vehicle's \ac{TCO} is composed of two main contributions: vehicle acquisition price $C_{\mathrm{a},i}$ and operation cost $C_{\mathrm{op},i}$
\begin{equation}\label{eq:TCOdiv}
	j_{\mathrm{TCO},i} = C_{\mathrm{a},i} + C_{\mathrm{op},i},
\end{equation}
where $C_{\mathrm{a},i}$ can be related to the production costs $C_{\mathrm{p},i}$ by means of a fixed overhead costs percentage $k_{\mathrm{oh}}$~\cite{UScommerce1998,Deutsprice}, as displayed in Fig.~\ref{fig:Cost_Breakdown}
\begin{equation}
	C_{\mathrm{a},i} = \frac{C_{\mathrm{p},i}}{k_{\mathrm{oh}}},
\end{equation}
while the operation costs are computed by multiplying the lifelong vehicle energy consumption $E_{\mathrm{v},i}$ and the energy cost $c_{\mathrm{e}}$,
\begin{equation}
	 C_{\mathrm{op},i} = E_{\mathrm{v},i} \cdot c_{\mathrm{e}}.
\end{equation}
In turn, $E_{\mathrm{v},i}$ is computed from the distance-specific energy consumption $F_{\mathrm{v},i}$ from the powertrain model in Section \ref{subsec:VECM} and the vehicle lifetime $d_{\mathrm{v,lt}}$ expressed as a distance
\begin{equation}
	E_{\mathrm{v},i} = F_{\mathrm{v},i} \cdot d_{\mathrm{v,lt}}.
\end{equation}
The vehicle production cost can be divided further into different contributions owing to the glider $C_{\mathrm{g}}$, motors $C_{\mathrm{m}}$, and battery pack $C_{\mathrm{b}}$
\begin{equation}
	C_{\mathrm{p},i} = C_{\mathrm{g},i}+C_{\mathrm{m},i}+C_{\mathrm{b},i}.
\end{equation}
\subsubsection{Glider Cost}
The glider comprises the body, chassis, low-voltage electrical components, exterior, and interior~\cite{BrookerGonderEtAl2015,DelDuceGauchEtAl2014}.
We also include the cost of the gearbox in the glider, as we assume it is not influenced by the production volumes~\cite{Barnes1998}. 
For the context of our analysis, and without loss of generality, we will consider the glider's cost a constant term depending exclusively on the vehicle class~\cite{SimeuKimEtAl2021}.
On the other hand, motor and battery costs are functions of the modules' size, multiplicity on the specific vehicle, and production volumes.

\subsubsection{Motor Cost} \label{subsubsec:mc}
We consider the motor module cost to scale linearly with its peak power $\overline{P}_{\mathrm{m}}$, following several studies~\cite{KoenigNicolettiEtAl2021,Cuenca1995,KochhanFuchsEtAl2014,FriesKerlerEtAl2017,RajashekaraMartin1995,VyasCuencaEtAl1998}.
Hence, depending on the number of modules the vehicle is equipped with, the total production cost of the motors becomes
\begin{equation}
	C_{\mathrm{m},i} = c_{\mathrm{m}} \cdot S_\mathrm{m} \cdot N_{\mathrm{m},i} \cdot \overline{P}_{\mathrm{m,o}},
\end{equation}
where $c_{\mathrm{m}}$ is the (peak) power-specific motor cost and is composed of two terms
\begin{equation}\label{cm}
	c_{\mathrm{m}} = c_{\mathrm{m},y} + \frac{\lambda_\mathrm{m,1}}{(N_{\mathrm{m},i} \cdot N_{\mathrm{v}} \cdot N-1)^{\lambda_\mathrm{m,2}}}.
\end{equation}  
The first term $c_{\mathrm{m},y}$ includes all the costs independent of volumes, like materials and process efficiency, indirectly depending on the market and technology level, which varies every year~\cite{Anderman2019}. 
The second term takes into account the effect of the \ac{EoS}. 
For this contribution, we use an asymptotic fit to reckon all the costs that can be amortized on a variable number of produced items $N_\mathrm{v}$. 
The coefficients $\lambda_\mathrm{m,1}$ and $\lambda_\mathrm{m,2}$ are subject to identification, as reported in \ref{app:CM}.
%In this way, for production volume approaching infinity, the module's cost converges to the raw material cost.

\subsubsection{Battery Cost}\label{subsubsec:bc}
The battery pack is one of the most expensive parts of an electric vehicle's powertrain~\cite{FriesKerlerEtAl2017}.
We assume its cost to linearly scale with the maximum battery module capacity
\begin{equation}
	C_{\mathrm{b},i} = c_{\mathrm{b}} \cdot  S_\mathrm{b} \cdot N_{\mathrm{b},i} \cdot \overline{E}_\mathrm{b,o},
\end{equation}
where $c_{\mathrm{b}}$ is the capacity-specific cost and it is composed of two terms, in a similar fashion to the motor model
\begin{equation}\label{cb}
	c_{\mathrm{b}} = c_{\mathrm{b},y} + \frac{\lambda_\mathrm{b,1}}{(N_{\mathrm{b},i} \cdot N_{\mathrm{v}} \cdot N-1)^{\lambda_\mathrm{b,2}}},
\end{equation}
where the coefficients $\lambda_\mathrm{b,1}$ and $\lambda_\mathrm{b,2}$ of the asymptotic term are identified from data available in the literature, as shown in \ref{app:CM}.

\subsection{Vehicle Powertrain Model}\label{subsec:VECM}
\begin{figure*}[t]
	\centering
	\includegraphics[width=1.46\columnwidth]{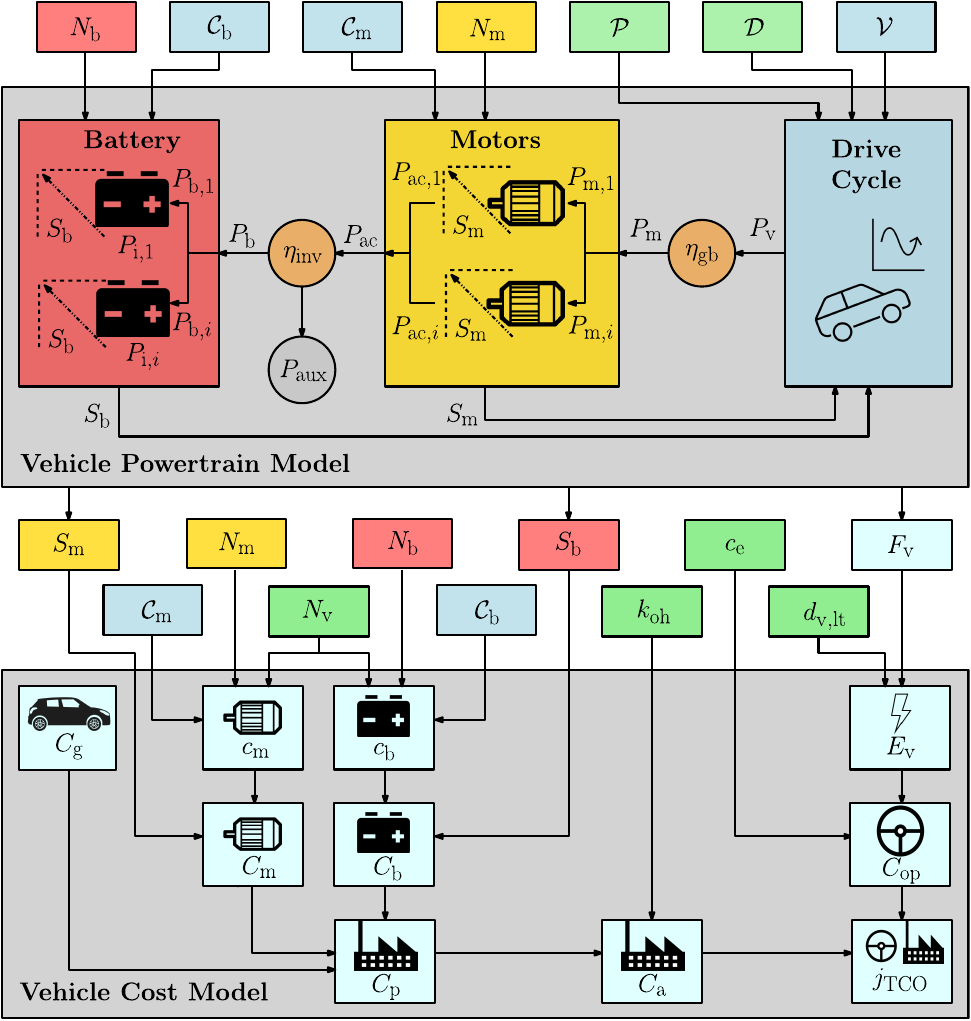}
	\caption{Flow diagram explaining the algorithm employed in computing the objective function (\ac{TCO}) of our concurrent design optimization framework.}
	\label{fig:model}
\end{figure*}

In this section, we describe the equations employed to model the powertrain behavior.
In accordance with common practices in this field~\cite{GuzzellaSciarretta2007,VerbruggenSalazarEtAl2019}, we apply a quasi-static modeling technique to predict the components' behavior and estimate the distance-specific energy consumption.
This formulation allows for the identification of the family-optimal modules' sizes for every configuration of modules in the family, accounting for their impact on energy consumption while still preserving convexity and its properties.
Fig. \ref{fig:model} displays the vehicle powertrain model as well as the cost model described in Section \ref{subsec:CM}, showing inputs, outputs, and interconnections.
Section~\ref{subsubsec:longitudinal} presents the vehicle's longitudinal dynamics, and Section~\ref{subsubsec:mass} gives insights into the mass model, taking into account the modules’ size and number.
Section~\ref{subsubsec:EM} introduces the EM model, whilst Section~\ref{subsubsec:BT} refers to the battery pack, and Section~\ref{subsubsec:perfconst} presents the required performance constraints that the vehicles need to satisfy.
For the sake of simplicity, we drop dependence on time $t$ whenever it is clear from the context.
Finally, even though we included extra passages in the derivation of some equations to facilitate understanding for the reader, we only label the version included in the formulation of \textbf{Problem \ref{prb:CDOP} (Concurrent Design Optimization Problem)} to avoid redundancy and maintain a clear mathematical definition of the problem.

\subsubsection{Longitudinal Dynamics}\label{subsubsec:longitudinal}
We optimize the design of the vehicles for a given driving cycle $\mathcal{D}$ of length  $d$ with exogenous longitudinal speed, acceleration, and slope trajectories $v(t)$, $a(t)$, and $\theta(t)$. %gathered in set $\mathbb{D}$ to lighten the notation.
Hence, the required power at the wheels $P_{\mathrm{v}}$ can be expressed as
\begin{multline}\label{eq:Preq}
	P_{\mathrm{v},i} = m_i \cdot v \cdot \left( c_{\mathrm{r},i} \cdot g \cdot \cos(\theta) + g \cdot\sin(\theta) + a \right) + \frac{1}{2} \cdot \rho \cdot c_{\mathrm{d},i} \cdot A_{\mathrm{f},i} \cdot v^3,
\end{multline}
where $m_i$ is the total mass of each vehicle, $g$ is the gravitational acceleration, $\rho$ is the density of the air, $c_{\mathrm{d},i}$ is the aerodynamic drag coefficient, and $c_{\mathrm{r},i}$ is the rolling resistance coefficient.

\subsubsection{Mass}\label{subsubsec:mass}
The vehicle's mass $m_{i}$ consists of several contributions: the glider mass $m_{\mathrm{g},i}$, the constant driver mass $m_{\mathrm{d}}$, the payload mass $m_{\mathrm{p},i}$, and the motor and battery mass, computed by scaling $m_{\mathrm{m,o}}$ and $m_{\mathrm{b,o}}$, accounting for the number of modules % characteristic of every vehicle type,
\begin{equation}\label{eq:mass}
	m_{i} = m_{\mathrm{g},i} + m_{\mathrm{d}} + m_{\mathrm{p},i} + m_{\mathrm{m.o}} \cdot S_{\mathrm{m}} \cdot N_{\mathrm{m},i} + m_{\mathrm{b,o}} \cdot S_{\mathrm{b}} \cdot N_{\mathrm{b},i}.
\end{equation}

\subsubsection{Transmission}\label{subsubsec:GB}
\begin{figure}[h]
	\centering
	\includegraphics[width=0.76\columnwidth]{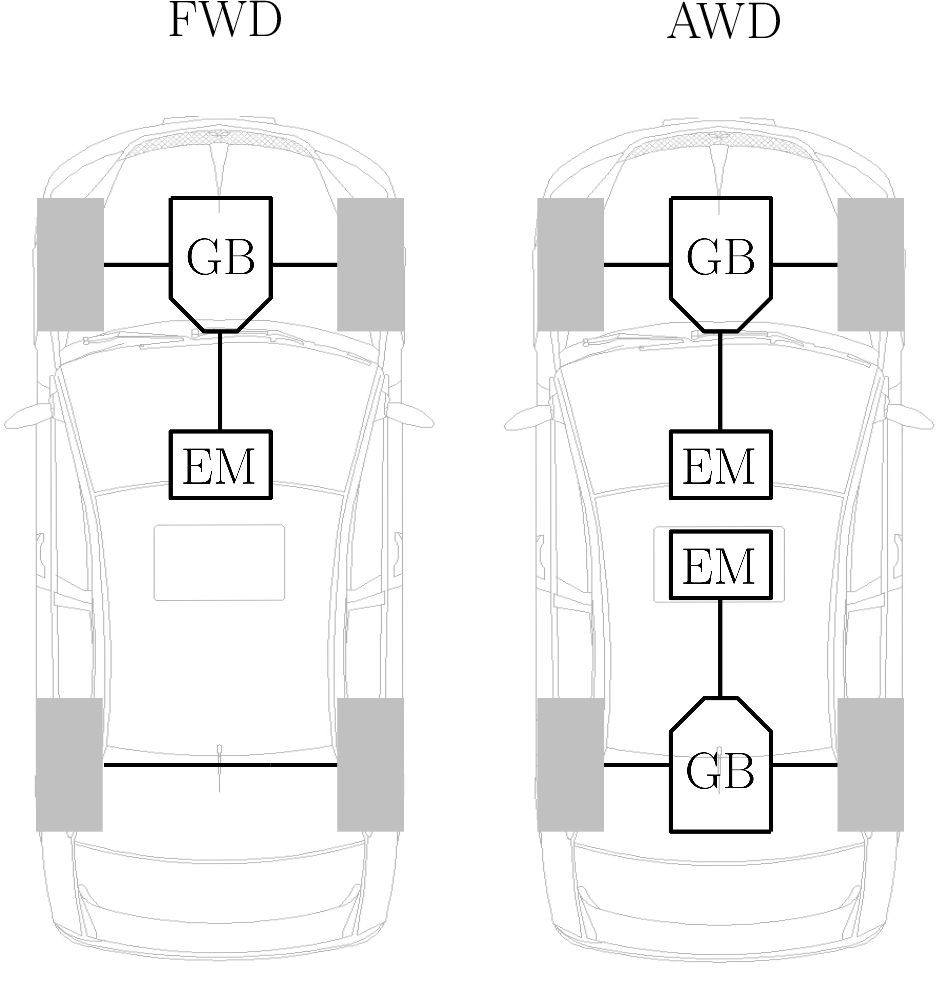}
	\caption{Vehicle topologies considered in this framework: Front-wheel drive (FWD) for vehicles equipping one motor module, and All-wheel drive (AWD) for two.}
	\label{fig:GB}
\end{figure}

We assume the regenerative braking system always in function when the power required at the wheels $P_{\mathrm{req}}$ is negative and it is equally divided among the $N_\mathrm{m}$ motor modules through fixed-gear transmissions (Fig. \ref{fig:GB}), whereby we consider the gear ratio $\gamma$ specifically designed for the motor speed interval considered.
The transmissions introduce losses that we model via a constant efficiency $\eta_\mathrm{gb}$:

%a number of fixed gear transmissions transfer to the $N_\mathrm{m}$ motor modules
%The power required at the wheels $P_{\mathrm{req}}$ is equally divided among the $N_\mathrm{m}$ motor modules through fixed-gear transmissions of ratio $\gamma$,

\begin{align*}
	P_{\mathrm{m}, i} = 
	\begin{cases}
		\frac{P_{\mathrm{v},i}}{\eta_\mathrm{gb} \cdot N_{\mathrm{m},i}} \quad & \text{if } P_{\mathrm{req},i}\geq 0 \\
		r_{\mathrm{b},i} \cdot \frac{\eta_\mathrm{gb} \cdot P_{\mathrm{v},i}}{N_{\mathrm{m},i}} \quad & \text{if } P_{\mathrm{req},i} < 0
	\end{cases},
\end{align*}
\begin{align*}
	r_{\mathrm{b},i} = \begin{cases} r_{\mathrm{FWD}} \quad & \text{if FWD}  \\
		r_{\mathrm{AWD}} \quad & \text{if  AWD} \end{cases},
\end{align*}
where the \rwmargin{regenerative braking fraction $r_\mathrm{b}$ accounts for the braking capabilities}{R5:2} of the vehicle's topology $r_{\mathrm{FWD}}$ and $r_{\mathrm{AWD}}$, included in the set $\mathcal{C_\mathrm{m}}$.
\rwmargin{This fraction only determines the maximum amount of power that can be recuperated based on the axles linked to a recuperation device, while the effective power regenerated depends on the efficiency of the motor(s) at that specific operating point.}{R5:2}
We can relax this constraint to an inequality without changing the problem solution (lossless relaxation) as
\begin{align}	\label{eq:Cdrive}
	P_{\mathrm{m},i} \geq \frac{P_{\mathrm{v},i}}{\eta_\mathrm{gb} \cdot N_{\mathrm{m},i}},\\
	P_{\mathrm{m},i} \geq r_{\mathrm{b},i} \cdot \frac{\eta_\mathrm{gb} \cdot P_{\mathrm{v},i}}{N_{\mathrm{m},i}}.
\end{align}
In fact, thanks to the particular problem structure, the optimal solution is located on the boundary, implying that these constraints will always hold with equality~\cite{EbbesenSalazarEtAl2018}.
Ultimately, assuming any value higher than the strictly necessary would be sub-optimal since it would entail higher losses and thus additional costs.

%The power is then  among one or more motors, depending on the topology considered
%\begin{equation}\label{eq:mass}
%	P_{\mathrm{m}} = \sum_{i=1}^{N}	P_{\mathrm{m}, i}.
%\end{equation}
\subsubsection{Electric Motor}\label{subsubsec:EM}
Analogously to our previous work~\cite{ClementeSalazarEtAl2022}, we make use of a second-order polynomial approximation of the reference motor losses $P_\mathrm{m,o,loss}$, extending a quadratic approximation~\cite{VerbruggenSalazarEtAl2019} to retain accuracy, without losing convexity
\begin{equation*}\label{eq:Ploss0}
	P_\mathrm{m,o,loss} = P_\mathrm{0}(\omega) + \beta(\omega) \cdot 	P_{\mathrm{m,o}} + \alpha(\omega) \cdot 	P_{\mathrm{m,o}}^2.
\end{equation*}
Furthermore, we assume that the operational limits $\underline{P}_{\mathrm{m,o}}$ and $\overline{P}_{\mathrm{m,o}}$ scale linearly with respect to the reference values
\begin{equation}\label{eq:Smrange}
	P_{\mathrm{m},i} \in \left[\underline{P}_{\mathrm{m,o}}, \overline{P}_{\mathrm{m,o}} \right] \cdot S_\mathrm{m}.
\end{equation}
Following the same rationale, we can write $P_{\mathrm{m,loss}}$ for every motor module as
\begin{equation*}
	P_{\mathrm{m,loss},i} = \left(P_\mathrm{0}(\omega) + \beta(\omega) \cdot \frac{P_{\mathrm{m},i}}{S_\mathrm{m}}  + \alpha(\omega) \cdot \frac{P_{\mathrm{m},i}^2}{S_\mathrm{m}^2} \right) \cdot S_\mathrm{m},
\end{equation*}
yielding
\begin{equation*}
	P_{\mathrm{m,loss},i} = P_\mathrm{0}(\omega) \cdot S_\mathrm{m} + \beta(\omega) \cdot P_{\mathrm{m},i} + \alpha(\omega) \cdot \frac{P_{\mathrm{m},i}^2}{S_\mathrm{m}}.
\end{equation*}
Finally, we obtain the AC input motor power $P_{\mathrm{ac},i}$ equation
\begin{equation*}
	P_{\mathrm{ac},i} = \left(P_{\mathrm{m},i}  + 	P_{\mathrm{m,loss},i}\right),
\end{equation*}
and we replace the expression for the losses
\begin{equation*}
	\label{eq:EMSOCC}
	P_{\mathrm{ac},i} = P_{\mathrm{m},i} +  P_\mathrm{0}(\omega) \cdot S_\mathrm{m} + \beta(\omega) \cdot P_{\mathrm{m},i} +  \alpha(\omega) \cdot \frac{P_{\mathrm{m},i}^2}{S_\mathrm{m}}.
\end{equation*}
Analogously to what has been done in Eq. \eqref{eq:Cdrive}, this equation can be losslessly relaxed for the purpose of retaining convexity~\cite{MurgovskiJohannessonEtAl2012}, obtaining the second-order conic constraint
\begin{equation*}
	\left(P_{\mathrm{ac},i} - P_{\mathrm{m},i} - S_\mathrm{m} \cdot  P_\mathrm{0}(\omega) - \beta(\omega) \cdot P_{\mathrm{m},i}  \right) + \frac{S_\mathrm{m}}{ \alpha(\omega)} \geq
\end{equation*}
\begin{equation} \label{eq:EMcone}
	\left \| \begin{matrix}
		2 \cdot P_{\mathrm{m},i} \\
		\left(P_{\mathrm{ac},i} - P_{\mathrm{m},i} - S_\mathrm{m} \cdot  P_\mathrm{0}(\omega) - \beta(\omega) \cdot P_{\mathrm{m},i} \right) - \frac{S_\mathrm{m}}{ \alpha(\omega)}		 
	\end{matrix}   \right \|_2.
\end{equation}

\subsubsection{Battery}\label{subsubsec:BT}
Similarly to what has been implemented for the EM, for each vehicle of the $i$-th type, every module's output battery power $P_{\mathrm{b},i}$ can be found starting from $P_{\mathrm{ac},i}$ by considering the inverter efficiency $\eta_\mathrm{inv}$, power consumption of auxiliary systems $P_{\mathrm{aux},i}$ (heating, air conditioning, lights, etc.), and motor and battery modules' multiplicity $N_{\mathrm{m},i}$ and $N_{\mathrm{b},i}$. Following the assumption that every module supplies an equal amount of output power in the $i$-th vehicle,
\begin{align*}
%	\label{eq:inverter}
	P_{\mathrm{b},i} = 
	\begin{cases}
		\frac{1}{N_{\mathrm{b},i}}\cdot \left(\frac{ P_{\mathrm{ac},i} \cdot N_{\mathrm{m},i}}{\eta_\mathrm{inv}} + P_{\mathrm{aux},i}\right) \quad & \text{if } P_{\mathrm{ac},i} \geq 0 \\
		\frac{1}{N_{\mathrm{b},i}} \cdot \left(\eta_\mathrm{inv} \cdot P_{\mathrm{ac},i} \cdot N_{\mathrm{m},i} + P_{\mathrm{aux},i}\right)\quad & \text{if } P_{\mathrm{ac},i} < 0
	\end{cases},
\end{align*}
that can be relaxed to
\begin{align}
	\label{eq:Cinverter}
	P_{\mathrm{b},i} \geq \frac{1}{N_{\mathrm{b},i}}\cdot \left(\frac{ P_{\mathrm{ac},i} \cdot N_{\mathrm{m},i}}{\eta_\mathrm{inv}} + P_{\mathrm{aux},i}\right), \\
	P_{\mathrm{b},i} \geq \frac{1}{N_{\mathrm{b},i}} \cdot \left(\eta_\mathrm{inv} \cdot 	P_{\mathrm{ac},i} \cdot N_{\mathrm{m},i} + P_{\mathrm{aux},i}\right).
\end{align}
We model the battery's internal losses using a function of the power that the battery would release when short-circuited: the ``short circuit power"~\cite{ClementeSalazarEtAl2022} $P_{\mathrm{sc},i}$.
We make use of a convex piece-wise affine approximation of the function composed of $K$ parts, depending on the energy $E_{\mathrm{b},i}$ and size $S_\mathrm{b}$~\cite{VerbruggenSalazarEtAl2019}.
%, derived in Appendix \ref{app:OCP}.
Hence, for every battery module,
\begin{equation*}
	P_{\mathrm{i},i} = P_{\mathrm{b},i} + \frac{P_{\mathrm{i},i}^2}{P_{\mathrm{sc},i}},
\end{equation*}
where
\begin{align*}%\label{eq:Psc}
	P_{\mathrm{sc},i} = \min_{k\in\{1,...,K\}} \left\{ a_k \cdot E_{\mathrm{b},i} + b_k \cdot S_\mathrm{b} \right\}.
\end{align*}
However, the form of this constraint would not allow a convex formulation, therefore, we apply a lossless relaxation to treat a set of affine inequalities.
\begin{equation}\label{eq:PscIN}
	P_{\mathrm{sc},i} \leq a_k \cdot E_{\mathrm{b},i} + b_k \cdot S_\mathrm{b}\;\forall k.\in\{1,...,K\}.
\end{equation}
Following the same logic applied in Section \ref{subsubsec:EM}, we also relax the battery losses equation to %linking $P_{\mathrm{b},i}$ to $P_{\mathrm{ac},i}$ and $P_{\mathrm{i},i}$ to $P_{\mathrm{b},i}$
\begin{equation}\label{eq:BTcone}
	\left(P_{\mathrm{i},i} - P_{\mathrm{b},i} \right) + P_{\mathrm{sc},i} \geq 
	\left \| \begin{matrix}
		2 \cdot P_{\mathrm{i},i} \\
		\left(P_{\mathrm{i},i} - P_{\mathrm{b},i} \right) - P_{\mathrm{sc},i} 		 
	\end{matrix}   \right \|_2.
\end{equation}
The overall energy of the battery $E_{\mathrm{tot},i}$ is determined by each module's size and the total number of modules per vehicle
\begin{equation}
	E_{\mathrm{tot},i} \in \left[ \ \overline{E}_\mathrm{b,o} \cdot \underline{\xi}, \overline{E}_\mathrm{b,o} \cdot \overline{\xi} \  \right] \cdot S_\mathrm{b} \cdot N_{\mathrm{b},i},
	\label{eq:Eblimits}
\end{equation}
while the dynamic of $E_{\mathrm{tot},i}$ is influenced by the internal power via
\begin{equation}\label{eq:diff}
	\frac{\mathrm{d} E_{\mathrm{tot},i}}{\mathrm{d} t} = - N_{\mathrm{b}} \cdot P_{\mathrm{i},i}.
\end{equation}
Ultimately, we evaluate the distance-specific energy consumption $F_{\mathrm{v},i}$ by dividing $E_{\mathrm{tot},i}$ by driving cycle length $d$
\begin{equation}
	F_{\mathrm{v},i}= \frac{E_{\mathrm{tot},i}(0) - E_{\mathrm{tot},i}(T)}{d}.
\end{equation}
Since we consider a variable efficiency of the battery through $P_{\mathrm{sc},i}$, we include a constraint to ensure that the operations are conducted around the half-capacity level of the battery by averaging the maximum capacity at the start and the minimum capacity at the end of the driving cycle
\begin{equation*}
	E_{\mathrm{tot},i}(0) = S_{\mathrm{b}} \cdot \overline{\xi}\cdot \overline{E}_\mathrm{b,o} \cdot N_{\mathrm{b},i},
\end{equation*}
\begin{equation*}
	E_{\mathrm{tot},i}(T) = S_{\mathrm{b}} \cdot \underline{\xi} \cdot \overline{E}_\mathrm{b} \cdot N_{\mathrm{b},i},
\end{equation*}
	\begin{equation}
	E_{\mathrm{tot},i}(0) + E_{\mathrm{tot},i}(T) = S_{\mathrm{b}} \cdot \left(\overline{\xi} + \underline{\xi} \right) \cdot \overline{E}_\mathrm{b,o} \cdot N_{\mathrm{b},i}.
	\label{eq:Ebrange}
\end{equation}
%For the sake of brevity, we refrain from proving that these relaxations are lossless, as the reason lies in the same principle.

\subsubsection{Performance Constraints}\label{subsubsec:perfconst}
Besides the equations modeling the powertrain behavior during vehicles' operations, we include constraints on the performance of each vehicle, derived from [1], to ensure that the design of every single vehicle meets all expectation set.
Hence, we write in convex form the acceleration time, top speed, power gradability, torque gradability, and range constraints as
\begin{equation}\label{eq:accT}
	N_{\mathrm{m},i} \cdot S_{\mathrm{m}} \cdot	t_{\mathrm{a},i} \leq \frac{\omega_{\mathrm{r}} \cdot r_{\mathrm{w},i}^2 \cdot m_{i}}{\overline{T}_{\mathrm{m,o}} \cdot \gamma^2} + \frac{m_{i}\cdot \left( v_{\mathrm{f}}^2 + \frac{\omega_{\mathrm{r}}^2}{\gamma^2} \cdot r_{\mathrm{w},i}^2  \right)}{2 \cdot \overline{P}_{\mathrm{m,o}}},
\end{equation}
\begin{equation}\label{eq:Vmax}
	N_{\mathrm{m},i} \cdot S_{\mathrm{m}} \cdot \overline{P}_{\mathrm{m,o}} \geq \frac{1}{2} \cdot \rho \cdot c_{\mathrm{d},i} \cdot A_{\mathrm{f},i} \cdot \left(v_{\mathrm{t},i}\right)^3,
\end{equation}
\begin{equation}\label{eq:Pgrad}
	N_{\mathrm{m},i} \cdot S_{\mathrm{m}} \cdot \overline{P}_{\mathrm{m},o} \geq m_{i} \cdot g \cdot v_{\mathrm{m},i} \cdot \sin(\theta_i),
\end{equation}
\begin{equation}\label{eq:Tgrad}
	N_{\mathrm{m},i} \cdot S_{\mathrm{m}} \cdot \overline{T}_{\mathrm{m,o}}\cdot \gamma \geq m_{i} \cdot g \cdot r_{\mathrm{w},i} \cdot \sin(\theta_i),
\end{equation}
\begin{equation}\label{eq:Range}
	E_{\mathrm{tot},i}(0) - E_{\mathrm{tot},i}(T) \leq N_{\mathrm{b},i} \cdot S_{\mathrm{b}} \cdot \left(\overline{\xi} - \underline{\xi} \right) \cdot \overline{E}_\mathrm{b,o} \cdot \frac{d}{d_{\mathrm{r},i}},
\end{equation}
where $\overline{T}_{\mathrm{m,o}}$ is the maximum reference torque
\begin{equation*}
	\overline{T}_{\mathrm{m,o}} = \frac{\overline{P}_{\mathrm{m,o}}}{\omega_{\mathrm{r}}}.
\end{equation*}
The values of the acceleration time $t_{\mathrm{a},i}$, top speed $v_{\mathrm{t},i}$, uphill speed $v_{\mathrm{m},i}$, slope $\theta_i$, and range $d_{\mathrm{r},i}$ are bounded by the required performance in the set $\mathcal{P}$
\begin{equation}\label{eq:accTB}
	t_{\mathrm{a},i} <= \overline{t}_{\mathrm{a},i},
\end{equation}
\begin{equation}\label{eq:VmaxB}
	v_{\mathrm{t},i} >= \underline{v}_{\mathrm{t},i},
\end{equation}
\begin{equation}\label{eq:Vmin}
	v_{\mathrm{m},i} >= \underline{v}_{\mathrm{m},i},
\end{equation}
\begin{equation}\label{eq:slope}
	\theta_i >= \underline{\theta}_i,
\end{equation}
\begin{equation}\label{eq:range}
	d_{\mathrm{r},i} >= \underline{d}_{\mathrm{r},i}.
\end{equation}

%We gather all these performance parameters in the set $\mathbb{P}$.
%\begin{equation*}
%	\mathbb{P} = \left\{ t_{\mathrm{acc}}, v_{\mathrm{f}}, v_{\mathrm{max}}, D_\mathrm{R}, v_{\mathrm{min}}, \theta_\mathrm{max} \right\}.
%\end{equation*}

\subsection{Optimization Problem Formulation}\label{subsec:OPF}

We formulate the concurrent design optimization problem with the objective of minimizing the TCO of the vehicle family by choosing the sizing of the shared modules and their multiplicity as follows.
\begin{problem}[\textbf{Concurrent Design Optimization Problem}]\label{prb:CDOP}%Optimal Shared Module Sizing Problem
Given a family of battery electric vehicles with modular powertrains, the TCO-optimal shared modules' sizes and multiplicity for the whole family are the solution of
	\begin{equation*}
		\begin{aligned}
			& \qquad \min_{S_\mathrm{m},S_\mathrm{b},N_{\mathrm{m},i},N_{\mathrm{b},i}} J_\mathrm{TCO} \\
			&\text{s.t. }\text{ Shared Constraints \eqref{eq:Sm}-\eqref{eq:Nb}}&\\
			& \qquad \text{Vehicle Constraints \eqref{eq:TCOdiv}-\eqref{eq:range}} & \forall k \quad \forall i\\
%			& \qquad N_{\mathrm{m},i} \in \mathbb{M}\\
%			& \qquad N_{\mathrm{b},i} \in \mathbb{B}.
%			& \qquad \text{Performance Constraints \eqref{eq:accT}-\eqref{eq:Range}} & \forall i \\
		\end{aligned}
	\end{equation*}	
\end{problem}
This problem can be solved with global optimality guarantees in a nested fashion:
For given $N_{\mathrm{m},i}$ and $N_{\mathrm{b},i}$, this problem can be framed as a second-order conic program and rapidly solved to global optimality with standard algorithms.
Thus, we analyze each of the $N_\mathrm{c}$ possible configurations 
\begin{equation*}
	N_\mathrm{c} = \left(	\lvert \mathbb{M} \rvert \cdot 	\lvert\mathbb{B}	\rvert \right)^{N},
\end{equation*}
leveraging the polynomial solving time of the convex approach, and identify the globally optimal solution through an exhaustive search.

%\begin{problem}[\textbf{Concurrent Design Optimization Problem}]\label{prb:CDOP}
%Given a family of battery electric vehicles with optimized shared modular powertrain, the optimal multiplicity configurations $N_{\mathrm{m}}$ and $N_{\mathrm{b}}$ of every vehicle in the family are the solution of 
%	\begin{equation*}
%	\begin{aligned}
%		& \min_{N_{\mathrm{m}},N_{\mathrm{b}}} J_{\mathrm{TCO}}(S_\mathrm{m}^\star,S_\mathrm{b}^\star)\\
%		\text{s.t. } 
%	\end{aligned}
%	\end{equation*}

\subsection{Discussion}\label{subsec:discussion}
In this paragraph, we present some clarification on the assumptions and limitations of our work.
First, we conservatively assumed that every vehicle has a fixed percentage of overhead costs, yet the methodology is still sound for every cost model incorporating the effect of \ac{EoS}.
Second, we scale the electric motor mass linearly as a function of the maximum power and the battery size only by acting on the number of cells in parallel, thus changing its energy without altering the battery voltage.
These scaling methods are in line with high-level modeling approaches and optimal sizing design problems~\cite{ArouaLhommeEtAl2023}.
%Moreover, if the size is between 50\% and 200\% of the reference, the approximations are quite accurate~\cite{GrunditzThiringer2018}.
Third, we assume that every vehicle is equipped with the same transmission modules.
The gear ratio $\gamma$ is considered to be specifically designed for the motor speed range, which remains constant while linearly scaling the motor in torque.
For a fair comparison, $\gamma$ is kept constant both in the individual vehicle-tailored design and in the concurrent design optimization.
This assumption is in accordance with the common design practice among electric vehicle manufacturers of using a standard reduction gear.
Furthermore, we focus our attention on the economical aspect of the design, since it is one of the major issues preventing a larger uptake of electric mobility~\cite{IEA2020}.
However, the powertrain sizing also influences the environmental aspect of BEVs.
In fact, batteries place a heavy burden on the environment: Their production process influences different impact categories such as acidification, eutrophication, human toxicity, eco-toxicity, and resource depletion.
Yet, the major contribution in the battery life-cycle emissions is not the production process, but, the electricity consumption during the battery's life inside the vehicle~\cite{SimonWeil2013,TagliaferriEvangelistiEtAl2016,EllingsenSinghEtAl2016}.
In our concurrent design optimization methodology, we sacrifice some energy efficiency in favor of cheaper vehicles, increasing the lifetime energy consumption (larger batteries, heavier vehicles) and possibly leading to more emissions when compared to (individual) vehicle-tailored design.
Nevertheless, since even in the case of a very carbon-intense electricity mix, BEVs still outperform conventional fossil-fuel-powered vehicles in terms of emissions, increasing the share of electric vehicles directly translates into a reduction of the life cycle emissions and environmental impact.
To compare the influence of each component's sizing on the overall environmental impact of the vehicle, an extended life cycle assessment would be needed, however, it is beyond the scope of our paper.
Finally, we ensure the problem convexity by adopting a convex objective function and a convex domain.
Every constraint included in \textbf{Problem \ref{prb:CDOP}} delineates a convex space, and, since the intersection of convex spaces is still convex, the overall domain is convex~\cite{Boyd2007}.
Since the only integer variables of the problem are the motor and battery multiplicities, the problem could be solved as a Mixed Integer Second Order Conic Program (MISOCP).
Nevertheless, we apply a nested approach where we solve a SOCP for every combination of the modules multiplicity in each vehicle, leveraging the fast solving time of convex programs to retain information on the other sub-optimal configurations, additionally obtaining sensitivity information.

\section{Results}\label{sec:results}
In this section, we showcase our methodology with a benchmark design problem considering the family design of four different Tesla models: Model~S, Model~3, Model~X and Model~Y. %concurrent design optimization
We compare the results of the concurrent strategy with the leading technique applied in the literature and on the market, consisting in optimizing the vehicle’s design specifically to minimize its own \ac{TCO}, whereby the components’ sizes are optimized individually for the single vehicle.
However, by comparing Table~\ref{tab:results} with Table~\ref{table:student}, we observe minor differences between the vehicle-tailored design and the actual commercial vehicle design. These differences can be ascribed to uncertainties in the models (we do not have access to the exact Tesla efficiency maps) and other company strategic decisions.

%a vehicle-tailored design optimization, whereby the modules' sizes are optimized individually for a single vehicle. 
%, to quantify the impact of the concurrent design strategy on performance and costs. %for a wide range of market conditions.
%Afterwards, we extend the study to different production volumes and energy cost, assessing the effectiveness of the methodology in .
The vehicle-tailored design is computed with the same framework considering only one vehicle at the time, with $N_\mathrm{b} = 1$ and $ N_\mathrm{m} = 2 $ (AWD topology).
In order to minimize the errors owing to market fluctuations, we adjusted the figures for inflation~\cite{Inflation} and converted the dollars into euros where necessary~\cite{Conversion}.
% referring to their value on the 25 of October 2022
In our analyses, we refer to market and technology from the year 2020.
Table~\ref{tab:costpar} gathers the parameters used for the cost model, Table~\ref{tab:ref} shows the reference parameters of the modules, whereas Table \ref{tab:Vpar} contains vehicle parameters, performance parameters, and vehicle type fractions.
%Cost function parameters
\begin{table}[t]
	\centering
	\caption{Cost model parameters~\cite{MEACP2023,CBS23,Eurostat2022}}
	\label{tab:costpar}
	\begin{tabular}{l c l}\toprule
		\textbf{Parameter}  &   \textbf{Value}  &   \textbf{Unit}
		\\ \midrule
		$c_{\mathrm{b},2020}$ & 79 & \unit{EUR/kWh} \\
		$\lambda_\mathrm{b,1}$ & 3911 & \unit{EUR/kWh} \\
		$\lambda_\mathrm{b,2}$ & 0.3278 & - \\ \midrule
		$c_{\mathrm{m},2020}$ & 2.2 & \unit{EUR/kWp} \\
		$\lambda_\mathrm{m,1}$ & 537 & \unit{EUR/kWp} \\
		$\lambda_\mathrm{m,2}$ & 0.4524 & - \\ \midrule
		$C_{\mathrm{e}}$ &  0.4 & \unit{EUR/kWh}\\
		$N_{\mathrm{v}}$ &  200000 & - \\
		$d_{\mathrm{v,lf}}$ & 200000 & \unit{km}\\
		$k_\mathrm{oh}$ & 0.5350 & - \\
		$C_\mathrm{g}$ & 14736 & \unit{EUR}\\
		\bottomrule
	\end{tabular}
\end{table}
%Parameters V

\begin{table}[t]
	\centering
	\caption{Parameters from sets $\mathcal{V}$,$\mathcal{P}$,$\mathcal{W}$~\cite{Tesla2023,EVDatabase2022}}
	\label{tab:Vpar}
	\begin{tabular}{l c c c c c}\toprule
		\textbf{Tesla}  &   \textbf{Model S}  &   \textbf{Model 3} & \textbf{Model X} & \textbf{Model Y} & \textbf{Unit}
		\\ \midrule
		$m_{\mathrm{g}}$ & 1105 & 1069  & 1328 & 1205 &\unit{kg}\\ 	
		$\eta_\mathrm{gb}$ & 0.98 & 0.98 & 0.98 & 0.98 & - \\
		$\eta_\mathrm{inv}$& 0.96 & 0.96 & 0.96 & 0.96 & - \\
		$A_{\mathrm{f}}$ & 2.34 & 2.2  &2.59& 2.66 & \unit{m^2} \\ 
		$c_{\mathrm{d}}$ &  0.24 & 0.23 & 0.24 & 0.23 & - \\ 
		$r_\mathrm{w}$  & 0.3518  & 0.3353 & 0.3759 & 0.3560&\unit{m} \\ 	
		$c_{\mathrm{r}}$ & 0.007 & 0.007 & 0.007 & 0.007 & - \\
		$m_{\mathrm{p}}$ & 0 & 0 & 500 & 500 & \unit{kg}\\
		$m_{\mathrm{d}}$ & 80 & 80 & 80 & 80 & \unit{kg}\\
		$P_{\mathrm{aux}}$ & 500 & 500& 500 & 500& \unit{W}\\
		%		Tesla Cybertruck & 2.73 & 0.39 - 0.48 & 1294  & 0.4140 & 0.007 \\ 
		\midrule
		$\overline{t}_{\mathrm{a}}$ &  3.3 & 6.1 & 3.9 & 6.9 & \unit{s}\\
		$\underline{v}_\mathrm{f}$ & 100 & 100 & 100 & 100 & \unit{km/h}\\ 
		$\underline{v}_{\mathrm{t}}$ & 261 & 225 & 262& 217& \unit{km/h}\\ 
		$\underline{v}_{\mathrm{m}}$ & 10 & 10& 10& 10& \unit{km/h}\\ 
		$\underline{\theta}$ &  25 &  25&  25&  25& \unit{\%}\\
		$\underline{d}_\mathrm{r}$ & 460 & 405 & 455 & 350 & \unit{km}\\
		\midrule
		$w$ & 0.25 & 0.25  & 0.25 & 0.25 & -\\
		\bottomrule
	\end{tabular}
\end{table}

\begin{figure}
	\centering
	\includegraphics[width=0.74\columnwidth]{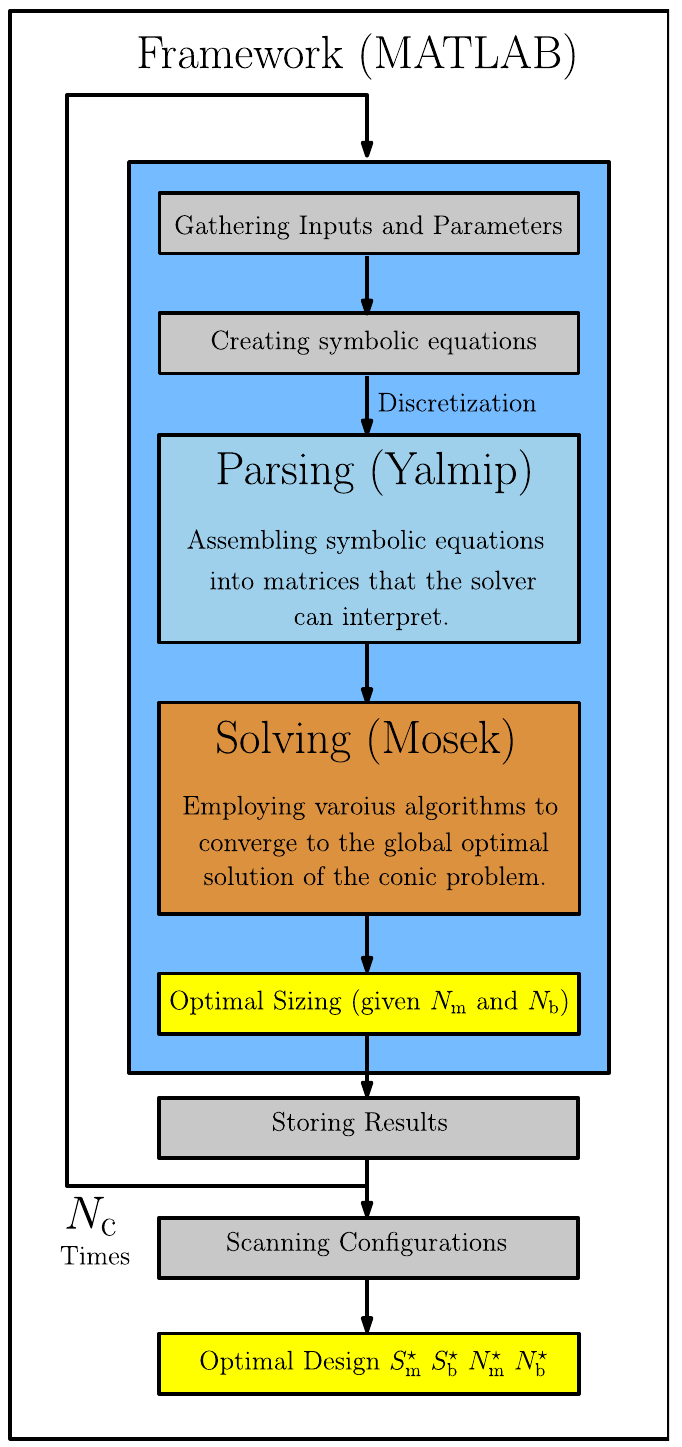}
	\caption{Flowchart of the solution process.}
	\label{fig:flowchart}
\end{figure}

We consider the Class 3 Worldwide harmonized Light-duty vehicles Test Procedure (WLTP) for the speed and acceleration trajectories and we discretize Problem \ref{prb:CDOP} using the forward Euler method with a sampling time of \mbox{1 \unit{s}}. Thereafter, we parse it with YALMIP~\cite{Loefberg2004} and solve it to global optimality with MOSEK~\cite{ApS2017}, in \mbox{2~\unit{s}} for each of the possible $N_\mathrm{c}$ combinations (Fig. \ref{fig:flowchart}).
%, for a total of approximately 7~\unit{h}.
It is important to underline that the solution depends on the sets $\mathbb{M}$ and $\mathbb{B}$.
In fact, they influence the possible ratios of motor power and battery capacity among the vehicles, thus potentially increasing efficiency and lowering costs.
Also, having more elements in the sets means a higher number of possible configurations, therefore increasing the total computation time.

Finally, in our benchmark problem, we consider equal fractions of vehicles of the $i$-th type, i.e., $ w_i =25\%$, giving the same importance to every vehicle type. 
\rwmargin{However, we perform a sensitivity analysis of the influence of this parameter on the vehicles' and family's TCO in Appendix C.}{R5:3}
%Since these parameters act like weights in the optimization, changing the fraction would inevitably prioritize the individual TCO reduction of the vehicle with higher weight at the expense of other types, as they would contribute in a larger share to the family TCO.
%If a vehicle's production volumes are much larger than other types', the benefit of sharing its modules with other vehicles to virtually increase the volumes is reduced.
%In the case of very different fractions, this benefit would be outweighed by the loss of energy efficiency and increased cost owing to oversized modules, required to satisfy the constraints of the other vehicles co-designed, causing a higher TCO compared to the vehicle-tailored design.

\begin{figure}[t]
						\centering
			\includegraphics[width=0.95\columnwidth]{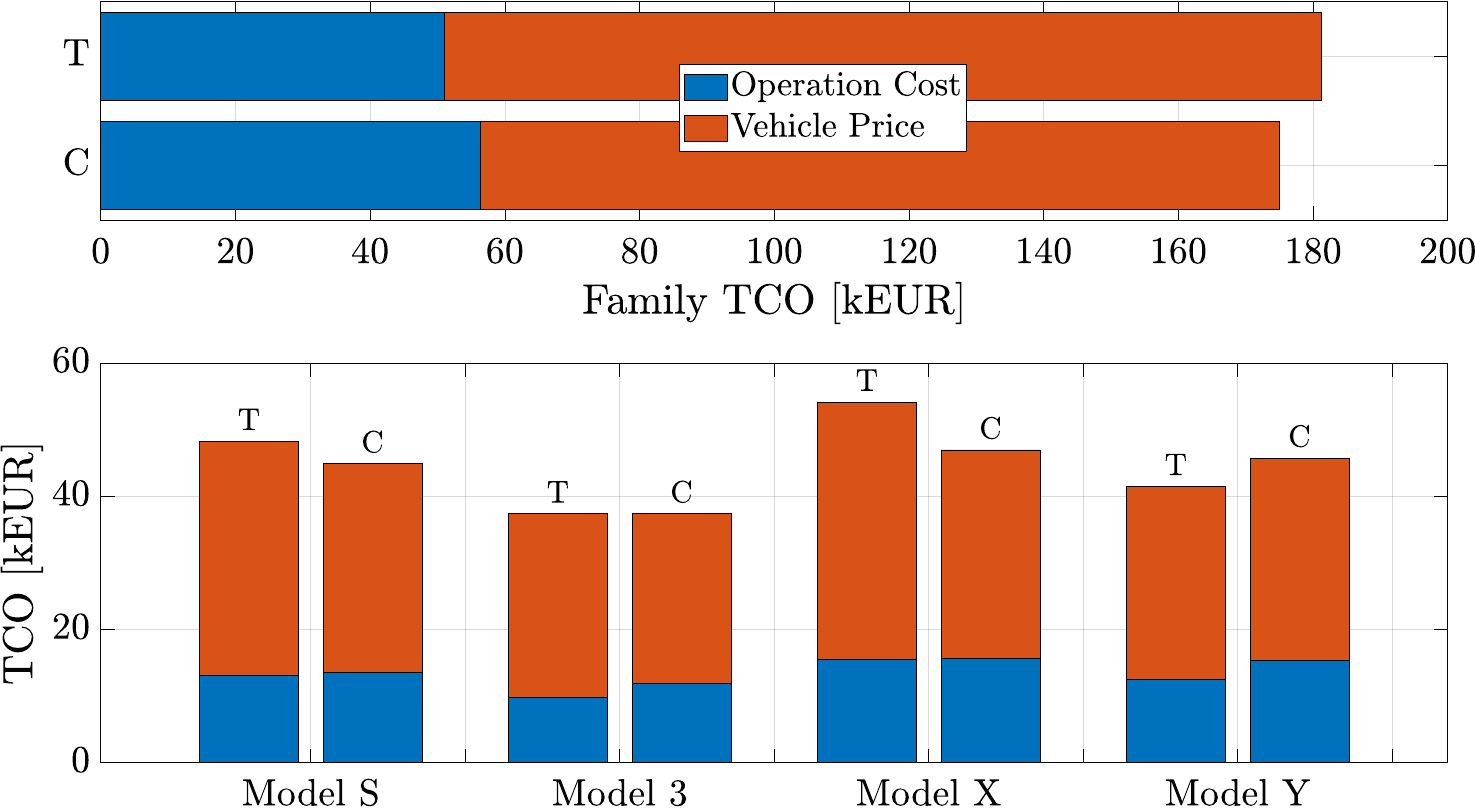}
			\caption{Comparison of costs between the concurrent design optimization (C) and the individual vehicle-tailored design (T). The bars on the top of the picture depict the overall \ac{TCO} of the family, showing a slight increment in operation cost against a significant reduction of the vehicle selling price. Below we show the single contribution of each vehicle.}
			\label{fig:cdesign}
		\end{figure}

Fig.~\ref{fig:cdesign} shows the comparison of the TCO achieved by concurrent and tailored design optimization for an energy price 40~EUR cents per kWh, and 200 thousand vehicles, highlighting that sizing the powertrain modules for the benchmark problem using a concurrent design optimization approach achieves a reduction of the Family \ac{TCO} of 3.5\% compared to the individual vehicle-tailored optimization. %, benefiting both manufacturers and customers.
%\begin{figure*}[t]
%	\begin{minipage}{0.49\linewidth}
%		\centering
%	\end{minipage}\hfill
%	\begin{minipage}{0.49\linewidth}
\begin{figure}
	\centering
	\includegraphics[width=\columnwidth]{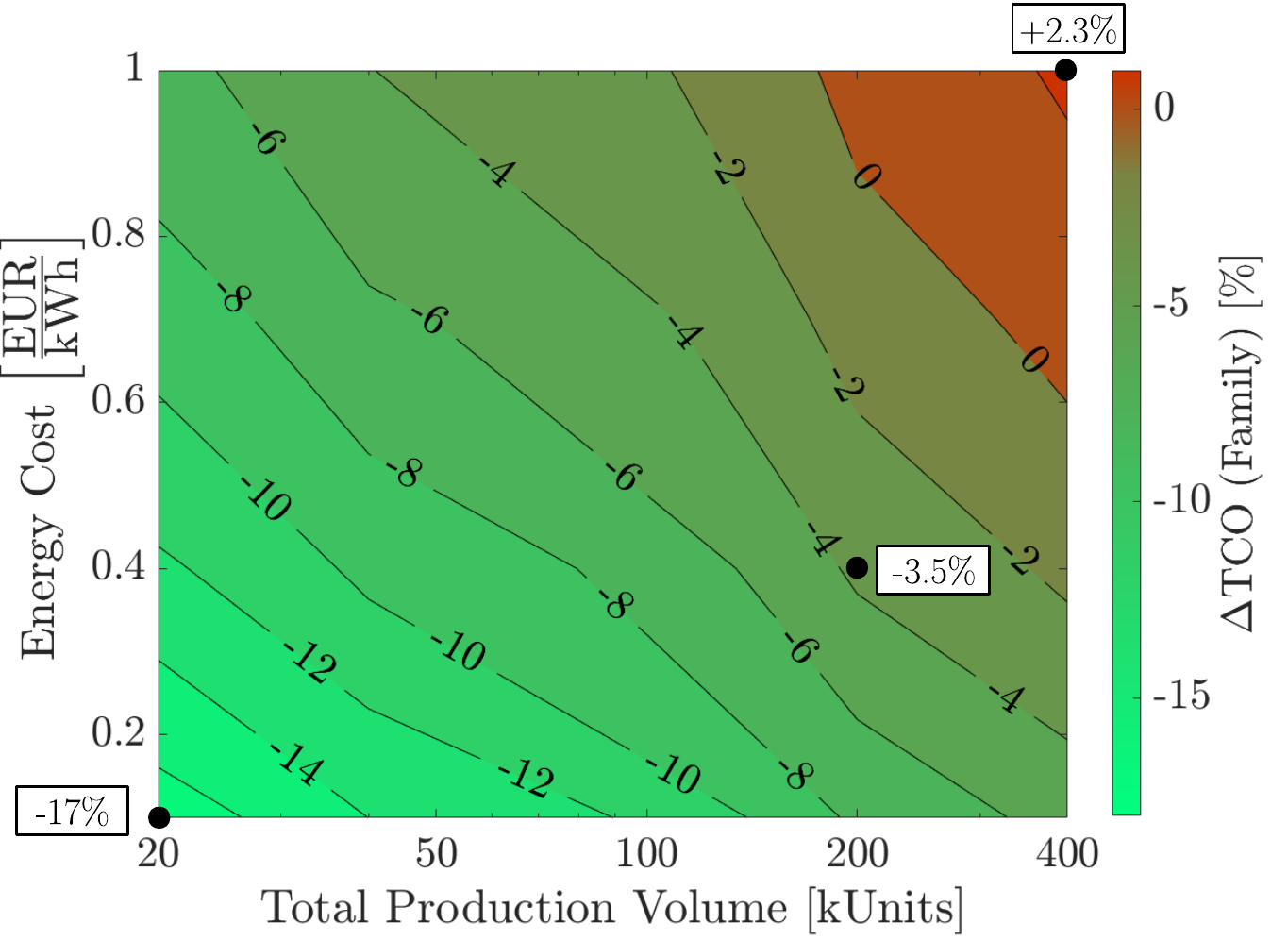}
	\caption{Profitability map of the concurrent design optimization approach for different market conditions of energy prices and production volumes. It is especially beneficial to share family-optimal modules for cheap energy prices and for small-scale production volumes.}
		\label{fig:scan}
			\end{figure}
%	\end{minipage}
%\end{figure*}
Specifically, there is an advantageous trade-off in using standardized modules to significantly reduce the cost of manufacturing while observing a limited increase in the operation costs thanks to the higher versatility prompted by modularity.
The vehicle price reduction is caused by the asymptotic relation linking the cost of an item and its production volume (\eqref{cm} and \eqref{cb}). For small volumes, the increase in the total amount of items caused by sharing the module features a significant item cost reduction. Conversely, for large volumes, the percentual cost reduction owing to the component sharing is smaller. For this reason, it may not keep pace with the increase in operation costs compared to the individual design.
In summary, although the traditional vehicle-tailored design converges to a more energy-efficient solution since the powertrain components are specifically optimized for that vehicle, the concurrent design optimization methodology produces a family design which is overall cheaper when considering the \ac{TCO}, thanks to the module cost reduction owing to the standardization and the \ac{EoS}.
When co-designing a large number of vehicles with very different performance constraints, some could present a higher TCO compared to the (individual) vehicle-tailored design. In fact, whilst always beneficial to the manufacturer, the family-optimal design could be disadvantageous to some users, ending up with a higher \ac{TCO} compared to its vehicle-tailored design if it leads to a greater cost reduction for the rest of the family (Model Y in Fig.~\ref{fig:cdesign}).
However, as long as the overall family TCO is lower, this strategy could be exploited to shift some of the costs from economy to luxury segments where a small increase in price would put less weight on the customers' choice.
Table~\ref{tab:results} reports the family design and performance specification of the concurrent design optimized family compared to its vehicle-tailored counterpart.
\begin{table}[t]
	\centering
	\caption{Benchmark problem family design specifications compared to a vehicle-tailored design}
	\label{tab:results}
	\begin{tabular}{l c  c  c  c c}\toprule
		\textbf{Tesla} & \textbf{Model S} & \textbf{Model 3}& \textbf{Model X} & \textbf{Model Y} & \textbf{Unit}
		\\ \midrule
		$j_\mathrm{TCO}$ & 44884 & 37381 & 46950 & 45730 & \unit{EUR}\\
		$\Delta$ & -6.98 & +0.04 & -13.28 & +10.22 & \unit{\%} \\
		\midrule
		$ m $ & 2232 & 1784 & 2455 & 2142 & \unit{kg}\\
		$\Delta $ & +6.57 & +7.57 & +0.60 & +16.31 & \unit{\%}\\
		\midrule
		$ d_\mathrm{r}  $ & 533 & 405 & 462 & 471 & \unit{km}\\
		$\Delta $ & +15.88 & 0.00 & +1.66 & +34.62 & \unit{\%}\\
		\midrule
		$t_\mathrm{a}$ & 3.28 & 4.96 & 3.90  & 6.38 & \unit{s}\\
		$\Delta$ & -0.71 & -18.72 & 0.00  & -7.60 & \unit{\%}\\
		\midrule
		$v_\mathrm{t}$ & 436 & 359 & 422 & 337 & \unit{km/h}\\
		$\Delta$ & +2.39 & +9.79 & +0.20 & +7.97 & \unit{\%}\\	
		\midrule
		$F_\mathrm{v}$ & 0.6101 & 0.5354 & 0.7031 & 0.6903 & \unit{MJ/km}\\
		$\Delta $ & +3.66 & +21.61 & +0.24 & +22.13 & \unit{\%}\\
		\midrule
		$\overline{P}_\mathrm{tot}$ & 625 & 313 & 625 & 313& \unit{kW} \\
		$\Delta$ & +7.33 & +32.33 & +0.60 & +25.87 & \unit{\%} \\
		\midrule
		$\overline{E}_\mathrm{tot}$ & 113 & 75 & 113 & 113&  \unit{kWh} \\
		$\Delta$ & +20.12 & +21.61 & +1.91 & +64.42 &  \unit{\%} \\
		\toprule 
		$S_\mathrm{m}$ &\multicolumn{4}{c}{2.33 (312 \unit{kW}) } & - \\
		$\Delta$ & +7.33 & +164.67 & +0.60 & +151.75 & \unit{\%}\\
		\midrule
		$S_\mathrm{b}$ &\multicolumn{4}{c}{1.60 (37.64 \unit{kWh})}& -  \\
		$\Delta$ & -59.96 & -39.20 & -66.03 & -45.19 & \unit{\%}\\
		\midrule
		$N_\mathrm{m}$ & 2 & 1 & 2 & 1 & - \\
		$\Delta$ & 0 & -50 & 0 & -50 & \unit{\%}\\
		\midrule
		$N_\mathrm{b}$ & 3 & 2 & 3 & 3 & - \\
		$\Delta$ & +200 & +100 & +200 & +100 & \unit{\%}\\
		\bottomrule
	\end{tabular}
\end{table}

We examine the sensitivity of the solution to different market conditions like energy price and production volume in Fig.~\ref{fig:scan}.
We consider electricity prices for household consumers in line with the statistics in the Netherlands~\cite{CBS23,MEACP2023} and other European countries~\cite{Eurostat2022}.
The electricity price linearly influences the operation cost, increasing the profitability of the concurrent design optimization methodology for affordable energy, whereby this strategy leverages a trade-off between energy efficiency and component costs. 
Furthermore, the relative profitability of the concurrent optimization methodology grows with decreasing production volumes up to a maximum theoretical limit, achieved for 1 unit per type produced.
Since the absolute saving in production cost grows asymptotically with production volumes (law of diminishing returns), the maximum theoretical limit does not have a practical application.
Therefore, in Fig.~\ref{fig:scan} we show typical production volumes of electric vehicle manufacturers from data of the latest years~\cite{Tesla2023b}.
Hence, on the one hand, employing a shared modular design is always advantageous for the vehicle manufacturer, who aims to produce as much as possible to leverage the cost reduction fostered by the \ac{EoS}.
On the other hand, the final user of the vehicle will benefit from a reduction in the acquisition price, which outweighs the increment in operation costs of a sub-optimal energy-efficiency family design.
Above a certain threshold, depending on the production volumes, this trade-off is no longer beneficial, and an individual vehicle-tailored design is more favorable.
Thus, the concurrent design methodology proves to be especially beneficial for small production volumes, becoming less and less advantageous for growing volumes, until the point where all the benefit gained from sharing the modules among the family is lost to compensate for the increase in cost of operation caused by the suboptimal energy-efficiency. Likewise, having cheap electricity prices favors the ``one size fits all” approach, lowering the significance of the cost of operation term.

%lower price of acquisition, which is balanced by an increased operation cost, depending on the electricity price.
%For very high energy costs the vehicle-tailored design is more favorable as the saving in price of acquisition 
%it is optimized for the energy efficiency of the specific vehicle

\begin{table}
	\centering
	\caption{Reference motor and battery parameters}
	\label{tab:ref}
	\begin{tabular}{l c c}\toprule
		\textbf{Symbol}  &   \textbf{Value}   &   \textbf{Unit}    
		\\	\midrule
		$\mathbb{M}$ &  $\left\{1, 2\right\}$ & -\\
		$\overline{S}_{\mathrm{m}}$  & 4 & - \\
		$\underline{S}_{\mathrm{m}}$  & 0.25 & - \\
		$\gamma$ &9.01 & - \\
		$m_{\mathrm{m,o}}$  & 81.6 & \unit{kg}\\ 
		$\overline{P}_{\mathrm{m,o}}$  & 134.08 & \unit{kW} \\
		$\omega_{\mathrm{r}}$ & 490.58 & \unit{rad/s}\\
		$r_{\mathrm{FWD}}$ & 0.6 & - \\
		$r_{\mathrm{AWD}}$ & 1 & - \\
		\midrule
		$\mathbb{B}$ & $\left\{1, 2, 3\right\}$ & - \\
		$\overline{S}_{\mathrm{b}}$  & 4 & - \\
		$\underline{S}_{\mathrm{b}}$  & 0.25 & - \\
		$m_{\mathrm{b,o}}$  & 138.6& \unit{kg}\\
		$\overline{E}_\mathrm{b,o}$  &23.48 & \unit{kWh} \\
		$\overline{\xi}$ & 0.9 & - \\
		$\underline{\xi}$ & 0.1 & -\\
		\bottomrule
	\end{tabular}
\end{table}

\section{Conclusions}\label{sec:conclusions}
This paper presented a concurrent design optimization framework to identify the optimal size and multiplicity of motor and battery modules that are installed within a family of electric vehicles, explicitly accounting for Economy of Scale (\ac{EoS}) effects on the resulting Total Cost of Ownership (\ac{TCO}). % \msmargin{assuming a fixed gearbox technology.}{too much}
To this end, we devised a cost model capturing the influence of production volumes, energy cost, modules' size and multiplicity on the vehicle's \ac{TCO}, and a powertrain model to account for its operations.
The resulting framework enables to compute minimum-\ac{TCO} vehicle family design solutions with global optimality guarantees.
%Our framework leverages a bi-level nested algorithm to jointly optimize the operation of the individual vehicles with the size and multiplicity of the electric motor and battery modules.
%The inner loop computes the optimal modules' size for every configuration, providing a solution that is guaranteed to be globally optimal with off-the-shelf second-order conic programming algorithms.
%The outer loop compares all the modules' configurations to identify the one featuring the lowest family \ac{TCO}.
When applied to a real-world case study for the design of the Tesla vehicle family, our methodology achieves a reduction of 3.5\% in the family \ac{TCO} compared to vehicle-tailored design solutions, whilst maintaining a comparable level of performance.
We showed that the profitability of the approach ranges from a reduction in family \ac{TCO} of 17\% for cheap energy prices and small production volumes to an increase of 2.3\% for very expensive energy prices and large production volumes.
%, where the increase in production brought by sharing modules minimally reduces the production costs, due to the asymptotic behavior of the cost function.
Even so, the approach is always beneficial to the vehicle manufacturer as it minimizes the production cost at the expense of the vehicles’ energy efficiency, while not requiring any particular cost investments.
% to a it delivers particular benefits (up to 17\%) for cheaper energy prices. %results especially beneficial
Overall, concurrent design optimization proved to be advantageous for a wide range of volumes and electricity prices, creating an opportunity to support the transition to electric mobility by substantially reducing the vehicle's upfront cost, one of the major factors hindering its adoption.
%In order to showcase our framework, we applied our concurrent optimization algorithm to the design of the Tesla vehicle family.
%Compared to the case where the modules are individually tailored to each vehicle, concurrently designing shared modules would reduce the family \ac{TCO} by almost 17\% while maintaining outstanding performance.
%Thereafter, we extended our analysis to different energy prices and production volumes, demonstrating the benefits of our design strategy for different market conditions, achieving improvements that range from X to Y\% in the case considered.
%In this context, on the one hand, vehicle manufacturers could exploit this methodology to set up the least amount of production lines possible while still enabling the design of a competitive vehicle family, able to meet different kinds of customer needs.
%On the other hand, end-users, would benefit from considerably lower vehicle acquisition prices while still preserving a good energy-consumption performance.
%Furthermore, sharing the same type of modules allows further advantages in logistics and availability of spare parts and would be especially advantageous for shared-vehicles fleet operators.
%, since it allows to keep using the modules in good condition from a vehicle at the end of its service life.

This work opens the field for the following extensions: First, we aim to sophisticate the framework by incorporating the optimization of the transmission and more powertrain architectures.
Second, we would like to assess the vehicles' environmental impact via life cycle analysis and study the effect of its minimization on the design solutions.
%Finally, AAA 
%include the gear ratio among the design variables. %Finally, it would be interesting to study the effect of using non-convex models and solvers on the solution.
\appendix

\section{Cost Model} \label{app:CM}

\begin{figure}
	\centering
	\includegraphics[width=\columnwidth]{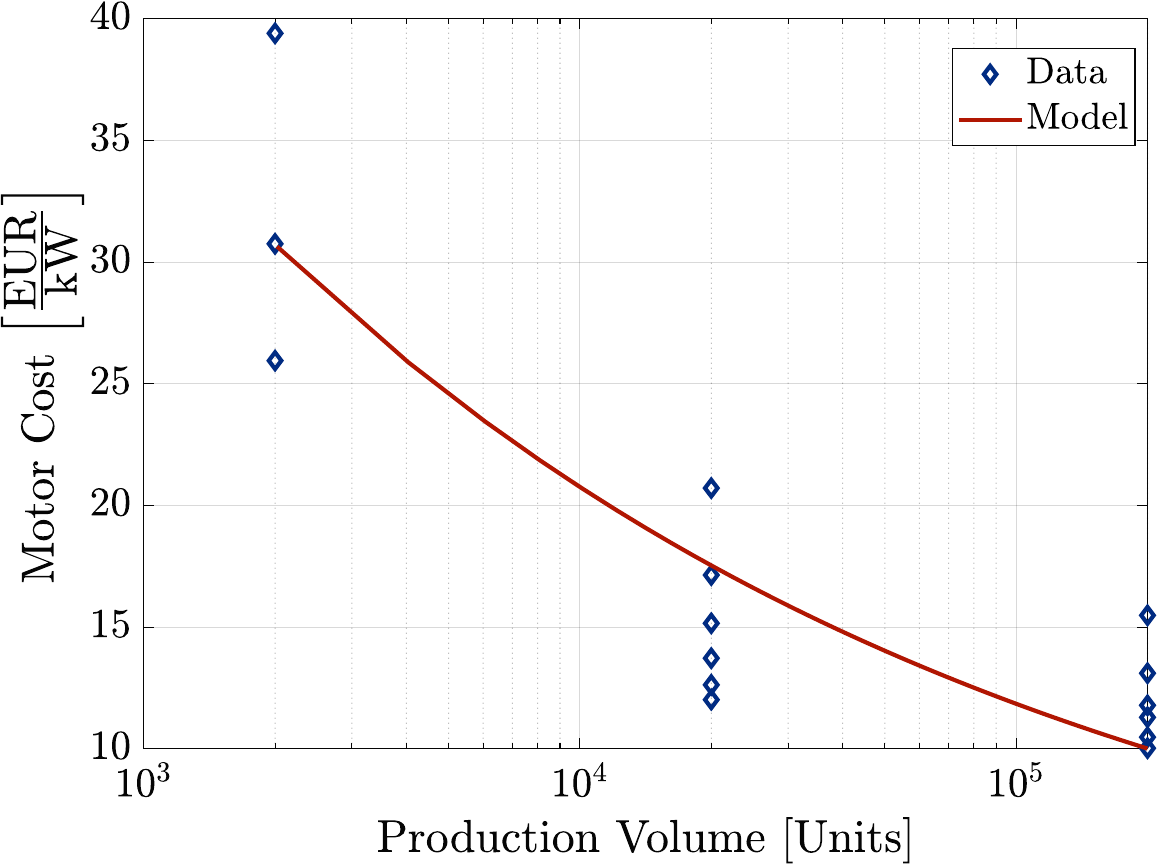}
	\caption{Production volumes impact on the motor cost in Euro per kW.}
	\label{fig:Masy}
	\end{figure}

\begin{figure}
		\centering
		\includegraphics[width=\columnwidth]{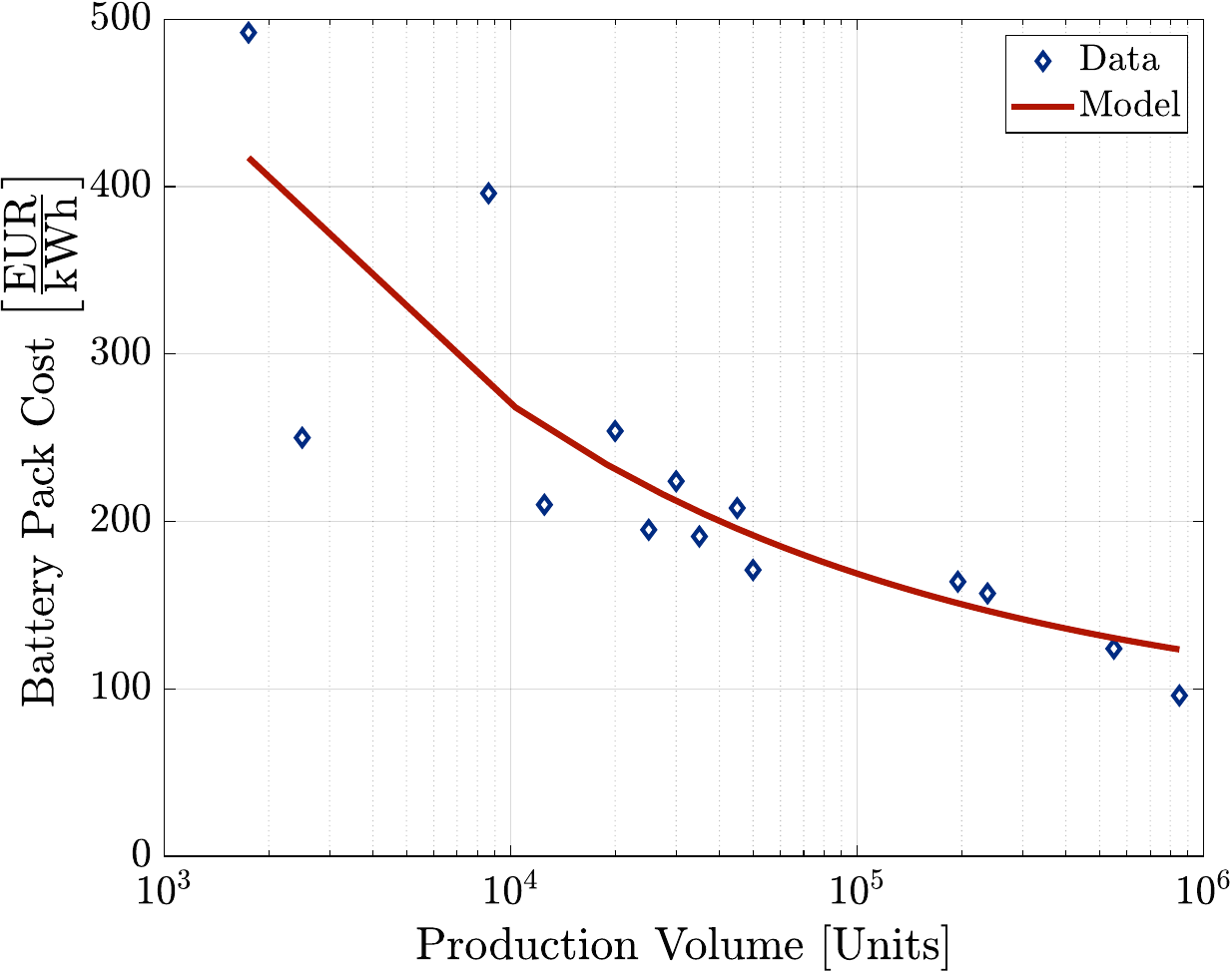}
		\caption{Production volumes impact on the battery pack cost in Euro per kWh.}
		\label{fig:Basy}
	\end{figure}

\begin{figure}
	\centering
	\includegraphics[width=0.95\columnwidth]{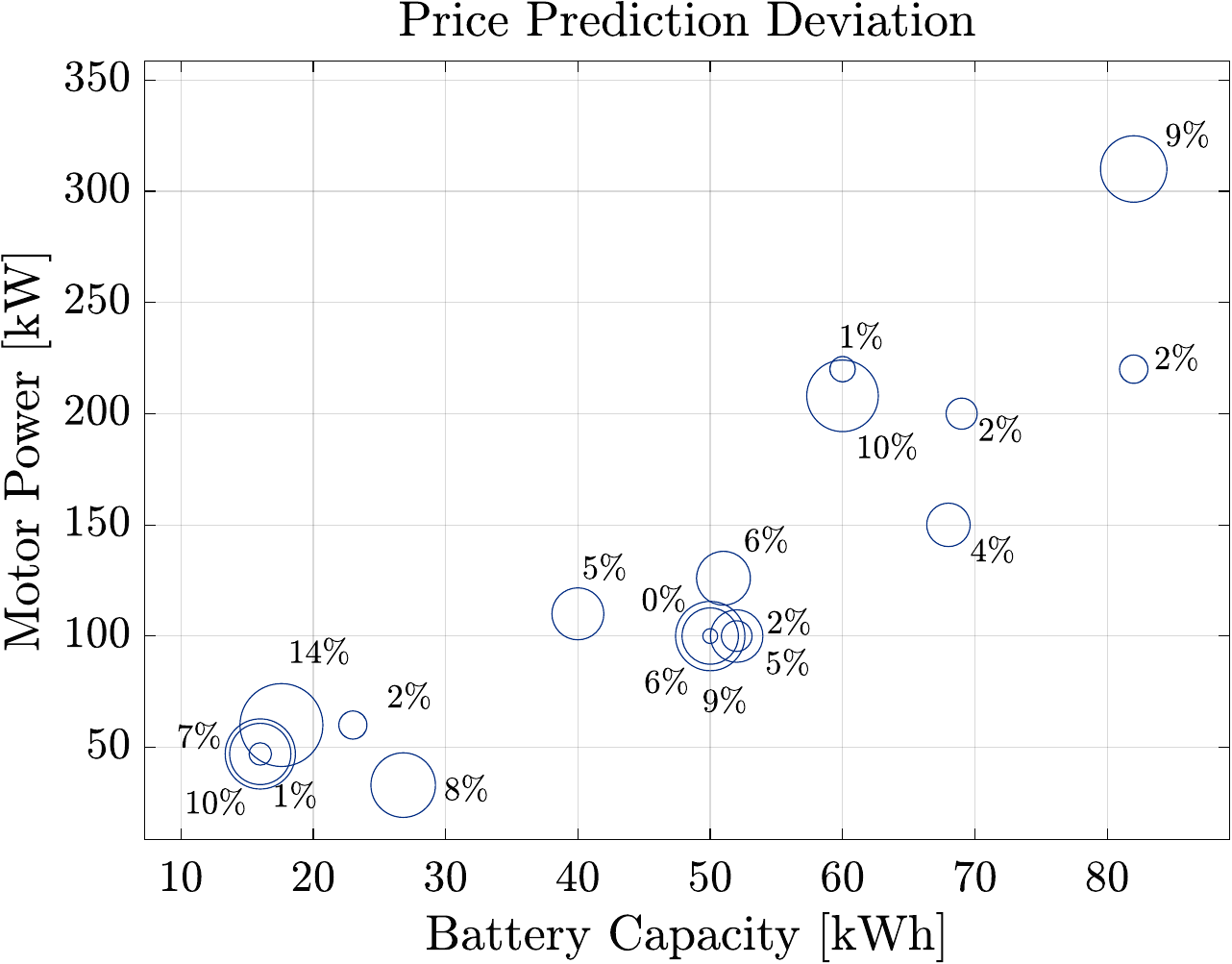}
	\caption{Percentage deviation of the model prediction from vehicles' price.}
	\label{fig:costcomparison}
\end{figure}

\begin{figure}
	\centering
	\includegraphics[width=0.8\columnwidth]{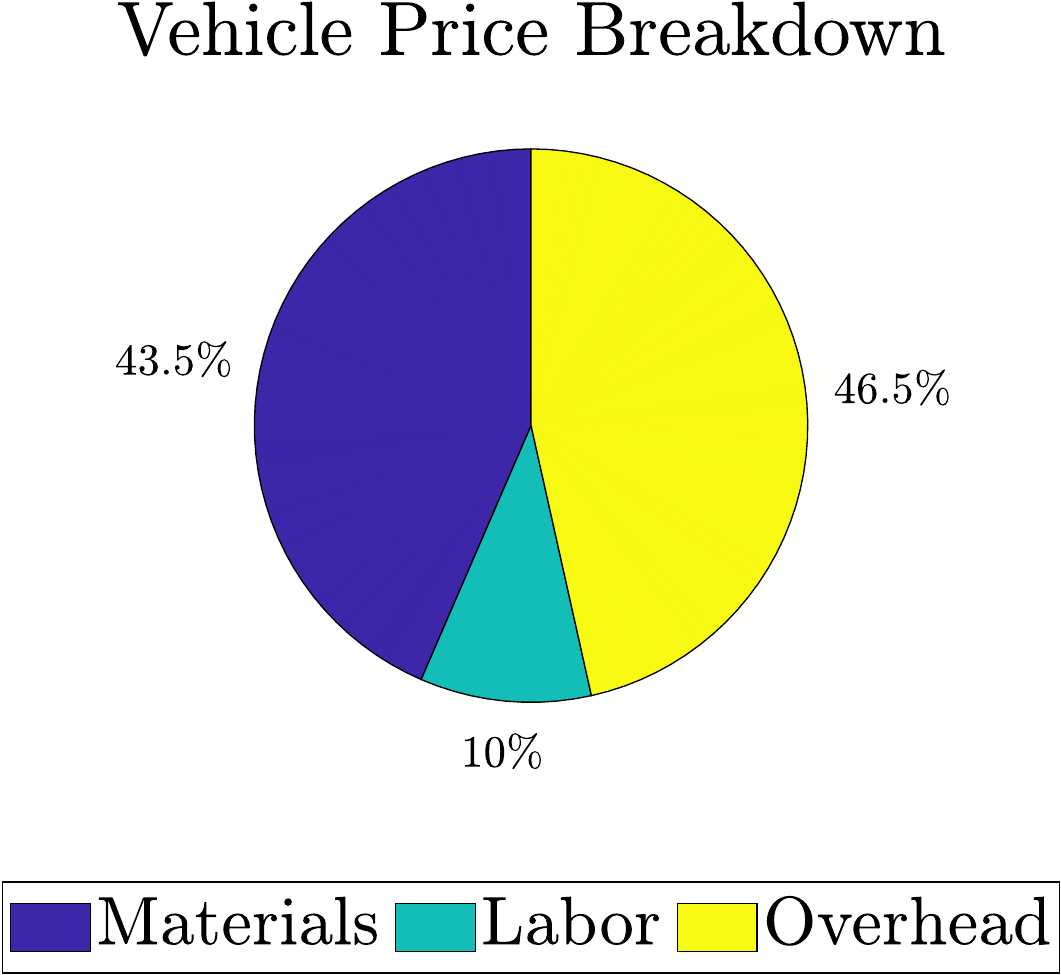}
	\caption{Price breakdown of a new vehicle \cite{KoenigNicolettiEtAl2021,Deutsprice}.}%,
	\label{fig:Cost_Breakdown}
\end{figure}

After a careful review of the literature, we developed a cost model able to portray the influence of \ac{EoS} strategies on battery packs and motor costs.
We report the data~\cite{FriesKerlerEtAl2017,KochhanFuchsEtAl2014,RajashekaraMartin1995,VyasCuencaEtAl1998,TechnologyAssessment1995,SantiagoBernhoffEtAl2012} used to identify $c_{\mathrm{b},2020}$, $\lambda_\mathrm{b,1}$, $\lambda_\mathrm{b,2}$, $c_{\mathrm{m},2020}$, $\lambda_\mathrm{m,1}$, and $\lambda_\mathrm{m,2}$ in Table~\ref{tab:data}, portraying the identified functions in Figures~\ref{fig:Masy} and \ref{fig:Basy}.
Furthermore, we indicate the glider cost depending on the vehicle type in Table~\ref{tab:glider_cost}.
Finally, we compute the selling price from the manufacturing cost by considering the overhead costs to be a constant fraction~\cite{KoenigNicolettiEtAl2021} $k_\mathrm{oh}$, as shown in Fig.~\ref{fig:Cost_Breakdown}.
The prediction of the cost model has been validated with a set of different \ac{BEVs}, achieving good results in modeling vehicle price using just the vehicle type, battery capacity, and motor power, assuming an average hundred thousand vehicles produced (Table \ref{tab:costcomparison} and Fig.~\ref{fig:costcomparison}). However, the error in the prediction increases for luxury brands, where the price is heavily influenced by factors which are not strongly coupled with manufacturing.

\begin{table}[b]
	\centering
	\caption{Glider Costs in 2020~\cite{SimeuKimEtAl2021,BEAN21}}
	\label{tab:glider_cost}
	\begin{tabular}{l c c}\toprule
		\textbf{Vehicle Type}  &   \textbf{Glider Cost [EUR]}    
		\\ \midrule
		City Cars & 7996 \\ %Compact
		Compact Cars & 10779 \\ %Mid-size
		%		 & 12718 \\ %Small SUV
		Large Cars & 14736 \\ %Mid-size SUV
		%		Pickup & 15964 \\
		\bottomrule
	\end{tabular}
\end{table}

%		\centering
%		\caption{Battery Pack capacity-specific cost by production volumes~\cite{KoenigNicolettiEtAl2021,Anderman2019,PhilippotAlvarezEtAl2019}}
%		\label{tab:Bdata}
%		\label{tab:Mdata}

%	\end{table}
%\end{minipage}\hfill
%\begin{minipage}{0.5\columnwidth}
%	\begin{table}[h]

%		\begin{tabular}{c c}\toprule
%			\textbf{Battery Pack}  &   \textbf{Cost} \\
%			\textbf{Units}  &   \textbf{[EUR/kWh]} \\ 
%			\midrule
%			1750 & 492  \\
%			2500 & 250 \\
%			8625 & 396 \\
%			12500 & 210 \\
%			20000 & 254\\
%			25000 & 195 \\
%			30000 & 224 \\
%			35000 & 191 \\
%			45000 & 208 \\
%			50000 & 171  \\
%			195000 & 164  \\			
%			237500 & 157 \\
%			550000 & 124 \\
%			850000 & 96 \\
%			\bottomrule
%		\end{tabular}
%	\end{minipage}

%\begin{figure*}[t]

%%	\end{minipage}
%\end{figure*}

\begin{table}
	\centering
	\caption{Peak-power-specific cost by production volumes~\cite{Lipman1999,Cuenca1995} (left), and Battery Pack capacity-specific cost by production volumes~\cite{KoenigNicolettiEtAl2021,Anderman2019,PhilippotAlvarezEtAl2019} (right)}
	\label{tab:data}
	%\begin{minipage}{0.5\columnwidth}
	\begin{tabular}{c c | c c}\toprule
		\textbf{Motor} & \textbf{Cost} &	\textbf{Battery Pack}  &   \textbf{Cost} \\
		\textbf{Units} & \textbf{[EUR/kW(p)]}  & \textbf{Units}  &   \textbf{[EUR/kWh]}\\
		\midrule  
		2000 & 39.40   & 1750 & 492  \\ %37
		2000 & 30.75   & 2500 & 250 \\
		2000 & 25.95   & 8625 & 396 \\
		20000 & 20.71  & 12500 & 210 \\
		20000 & 17.14  & 20000 & 254\\
		20000 & 15.15  & 25000 & 195 \\
		20000 & 13.72  & 30000 & 224 \\
		20000 & 12.62  & 35000 & 191 \\
		20000 & 12.01  & 45000 & 208 \\
		200000 & 15.48 & 50000 & 171  \\
		200000 & 13.10 & 95000 & 164  \\
		200000 & 11.78 & 237500 & 157 \\
		200000 & 11.29 & 550000 & 124 \\
		200000 & 10.47 & 850000 & 96 \\
		200000 & 10.02 & & \\ %55
		\bottomrule
	\end{tabular}
\end{table}

\begin{small}
	\begin{table}
		%	\begin{minipage}{0.49\linewidth}
			\caption{Validation of the vehicle's cost model assuming an average production volume of a hundred thousand vehicles per type}
			\label{table:student}
			\centering
			\resizebox{\columnwidth}{!}{%
			\begin{tabular}{l c c c c c}\toprule \label{tab:costcomparison}
				\textbf{Vehicle}  &   \textbf{Capacity}  &   \textbf{Power} & \textbf{Pred.} & \textbf{Price}\\
				&   \unit{\textbf{kWh}}  &  \unit{\textbf{kW}} & \unit{\textbf{EUR}} & \unit{\textbf{EUR}} & 
				\\ \midrule % Compact - CLASS AB
				\textbf{City (Seg. A-B)}\\
				Mitsubishi i-MiEV & 16  & 47& 20188 & 19990~\cite{I-MIEV2023} \\ 
				Peugeot iOn & 16 & 47 & 20188 & 22360~\cite{Database2023} \\
				Citroen C-Zero & 16 & 47& 20188 & 21800~\cite{Database2023}\\ 
				Smart EQ fortwo & 17.6 & 60& 20746 & 23995~\cite{Database2023} \\
				Renault Twingo E &  23  & 60 & 22450 & 22105~\cite{Twingo2023}  \\
				Dacia Spring El. & 26.8  & 33 & 23537 &  21750~\cite{DACIA2023}\\
				%		Abarth 500e Scoprionissima &42 kWh& 113 kW & 28663 & 28318 \\
				Opel Corsa-e & 50 & 100 & 31134 & 30999~\cite{OPEL2023} \\
				Peugeot e-208 & 50 & 100& 31134 &  33220~\cite{Peugeot2023} \\
				Renault Zoe Q90 & 52 & 100 & 31765 & 33590~\cite{Database2023}\\
				\midrule
				\textbf{Compact (Seg. C)}\\
				%		Fiat 500e Hatchback & 23.8 kWh & 70 kW & 28008 & 28990 & \cite{Pricelist500}\\
				Nissan Leaf & 40 & 110& 33222 & 35090~\cite{Database2023} \\
				%		Fiat 500e Hatchback & 42 kWh & 87 kW & 33758 & 34490 & \cite{Pricelist500} \\ 
				Citroen e-C4 X & 50  & 100 & 36336 & 40140~\cite{CITROEN2023}  \\
				ORA Funky Cat & 51  & 126 & 36758 & 38990~\cite{Database2023} \\
				Renault Zoe R135 & 52  & 100 & 36967 & 36295~\cite{Database2023} \\
				%		Opel Astra & 54 kWh &115 kW & 32458 & 36650 & \\
				%			Renault Zoe ZE50 R110 & 54.7 & 80 & 37736 & 34985~\cite{Database2023} \\		
				Kia Niro EV & 68&  150& 42221 & 43850~\cite{KiaNiro2023}\\
				\midrule
				\textbf{Large (Seg. D-E)}\\
				Tesla Model 3 & 60 & 208 & 47333 & 42993~\cite{Tesla2023} \\
				Tesla Model Y  & 60 &  220 & 47382 & 47993~\cite{Tesla2023} \\
				Polestar 2 & 69 & 200 & 50222 & 51200~\cite{Polestar2023} \\
				Polestar 2 LR & 82 & 220 & 54324 & 55200~\cite{Polestar2023}\\
				Polestar 2 Dual LR & 82 & 310 & 54695 & 59900~\cite{Polestar2023}\\
				%			BMW iX xDrive40 &  76.6 &  240 & 52703 & 54720\\
				%			BMW i4 eDrive40 &  82 &  &
				%https://www.polestar.com/nl/polestar-2/specifications/
				%			Tesla Model S  & 100 &  504 & 61175 & 94990~\cite{Tesla2023} \\
				%			Tesla Model X  & 100 &  504 & 61175 & 99990~\cite{Tesla2023} \\	 
				\bottomrule
				%		Mini Cooper SE & 32.6 kWh & 135 kW & 25788 & 38691   \\
				%		\\ \midrule % Exec - CLASS E
				%		Kia EV6GT & 74 kWh & 430 kW & 54962 & 57845 \\
				%		Tesla Model Y LongRng Dual  & 82kWh & 378kW & 57271 & 59111 \\
				%		BMW iX1 xDrive30 & 68.0 kWh&  230 kW & 52243 & 56568 \\ 
				%		BMW iX xDrive50 & 111.5 & kWh 385 kW & 66608 & 103755 \\ 
				%		BMW iX M60 & 111.5 kWh&  455 kW & 66896 & 135699 & \\ 
				%		Audi Q8 e-tron 55 quattro &   114.0 kWh 77 & 300 kW & 54604	& 85300 \\ 	 	    
				%		BMW i4 eDrive35 &70.2 kWh & 210 kW & 43305 & 58080 \\
				%		Honda e Advance & 35.5 kWh & 113 kW & 26426 & 40615 & 39900 \\
			\end{tabular}
		}
	\end{table}
\end{small}

\vspace{50pt}

\section{Powertrain Component Models} \label{app:PCM}
\begin{figure}
	\centering
	\includegraphics[width=0.9\columnwidth]{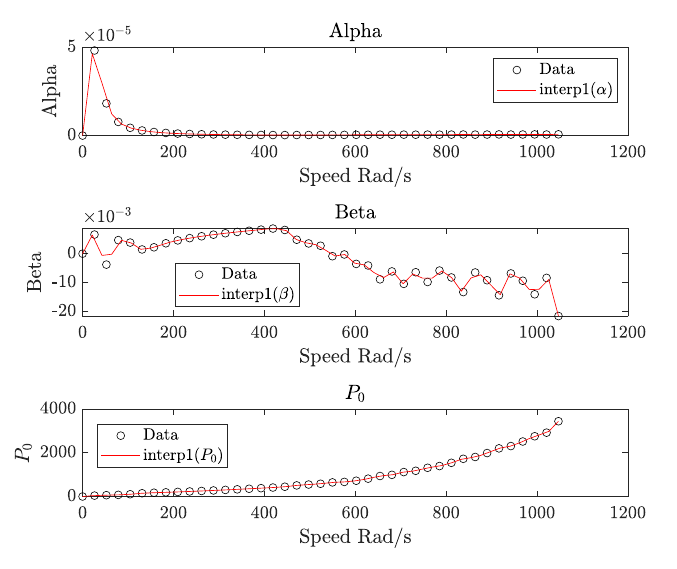}
	\caption{Speed-level dependent identified motor loss coefficients.}
	\label{fig:motor_coeff}
\end{figure}

\begin{figure}
	\centering
	\includegraphics[width=0.9\columnwidth]{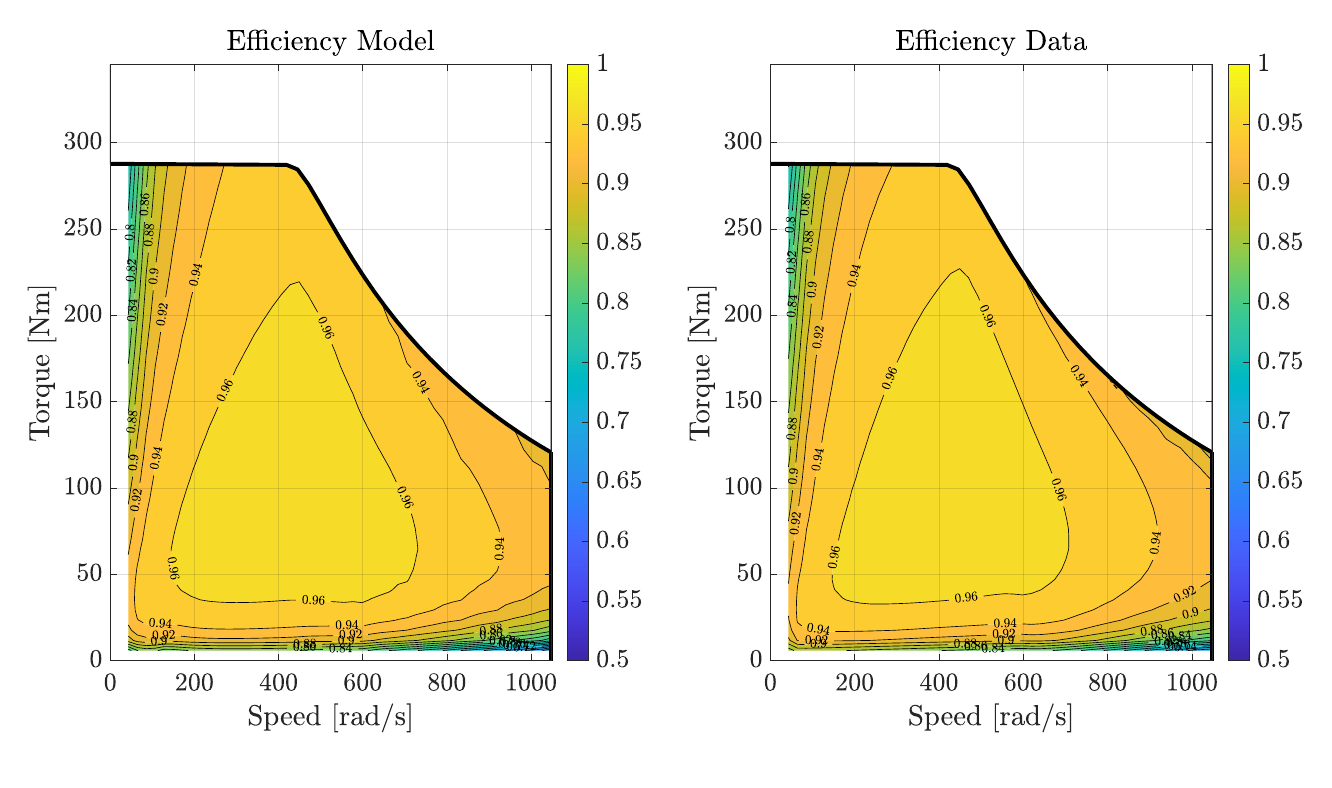}
	\caption{Reference electric motor efficiency map from the model (left) compared with data (right). Data from \cite{Solipuram2024}}
	\label{fig:motor_eff}
\end{figure}

\begin{figure}
	\centering
	\includegraphics[width=0.9\columnwidth]{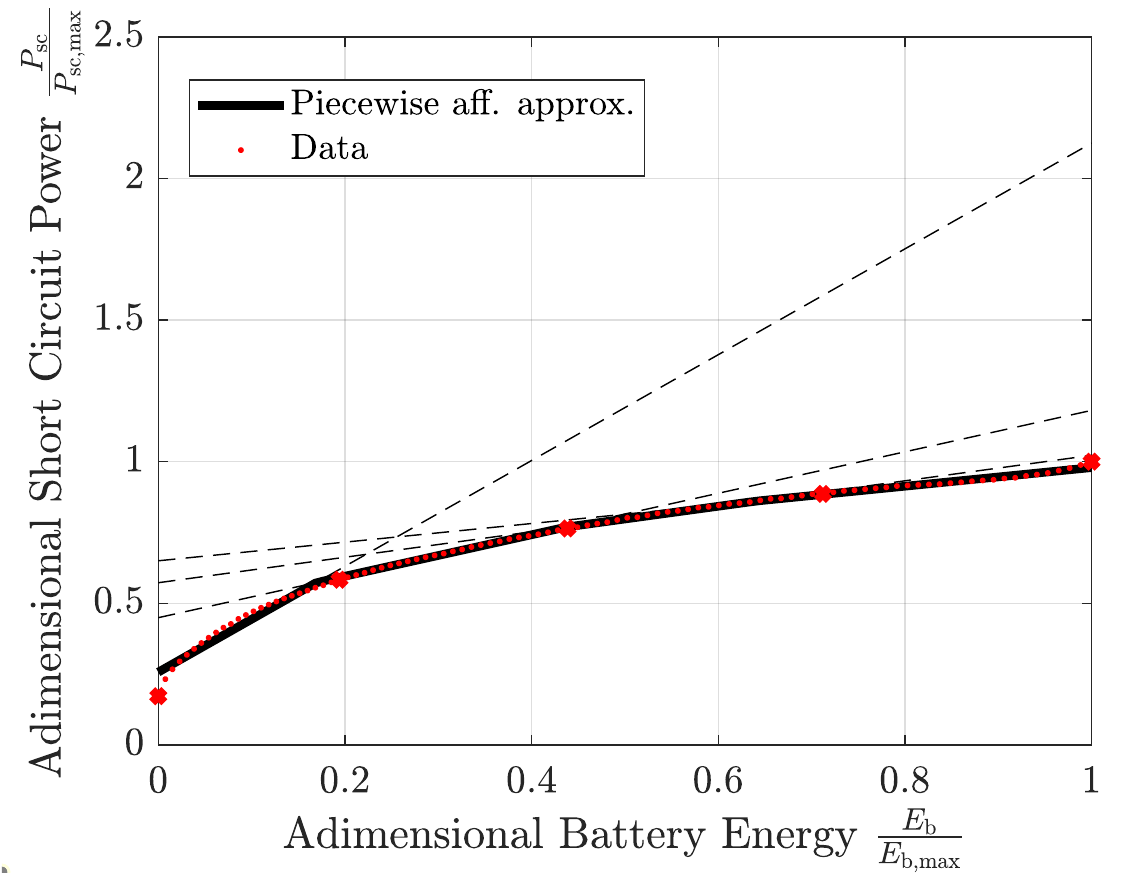}
	\caption{Reference battery short circuit power as a function of the battery energy. The slope and the intercept of the piecewise affine approximations are the battery loss coefficients $a_k$ and $b_k$.}
	\label{fig:battery_psc}
\end{figure}

The motor and the battery models capture the susceptibility of the efficiency to the components' sizing and the different dynamic driving conditions determined by the cycle $\mathcal{D}$, according to the equations presented in Sections \ref{subsubsec:EM} (Motor) and \ref{subsubsec:BT} (Battery).
We show the identified motor parameters of Eq.~\eqref{eq:EMcone}, namely $P_\mathrm{0}(\omega)$, $\beta(\omega)$ and $\alpha(\omega)$, in Fig.~\ref{fig:motor_coeff}, resulting in the efficiency map from Fig.~\ref{fig:motor_eff}, with a normalized root mean square error (NRMSE) of 0.41\% on the efficiency.
Likewise, the battery losses depend on the function $P_{\mathrm{sc}}$ presented in Eq.~\eqref{eq:BTcone}, which is also determined through a parameter identification process on the reference component, with an NRMSE of 1.32\% on $P_{\mathrm{sc}}$.
The slope and the intercept of the piecewise affine approximations of $P_{\mathrm{sc}}$ presented in Fig.~\ref{fig:battery_psc} are the battery loss coefficients $a_k$ and $b_k$, respectively, from Eq.~\eqref{eq:PscIN}.
Finally, in line with common practices in the field~\cite{GuzzellaSciarretta2007}, we assumed the gearbox and inverter efficiencies constant since the impact of different operative conditions on their efficiency is lower compared to other powertrain components, such as motor and battery.

\section{Vehicle Type Fraction Sensitivity Analysis} \label{app:VTFSA}

\begin{figure*}
	\centering
	\includegraphics[width=\linewidth]{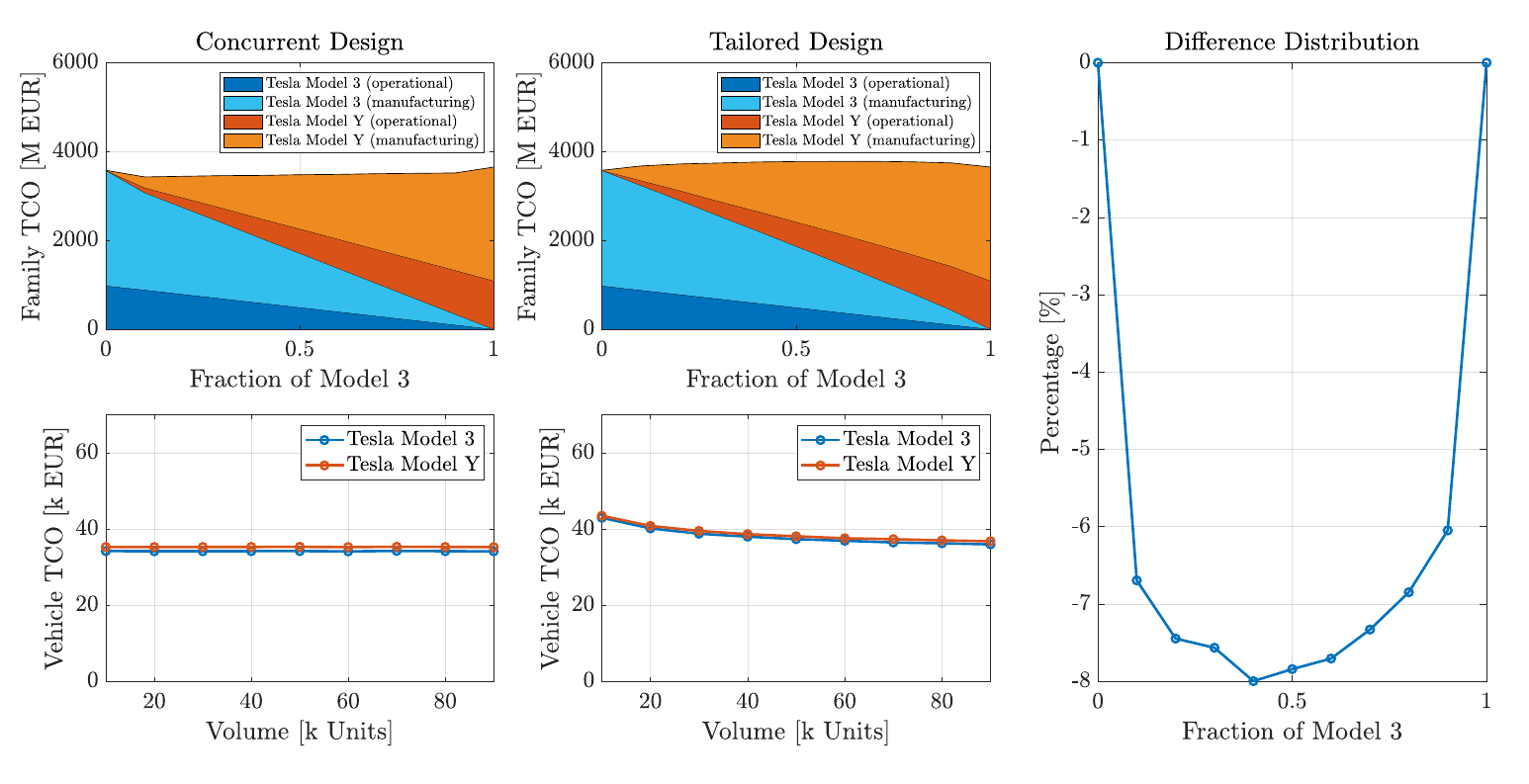}
	\caption{Sensitivity analysis of Tesla Model 3 and Tesla Model Y with respect to the vehicle type fraction $w_i$ and a total production of 100000 units. For fractions 0 or 1, one of the vehicles is not produced while the other is optimized individually.}
	\label{fig:ratio}
\end{figure*}

\begin{figure*}
	\centering
	\includegraphics[width=\linewidth]{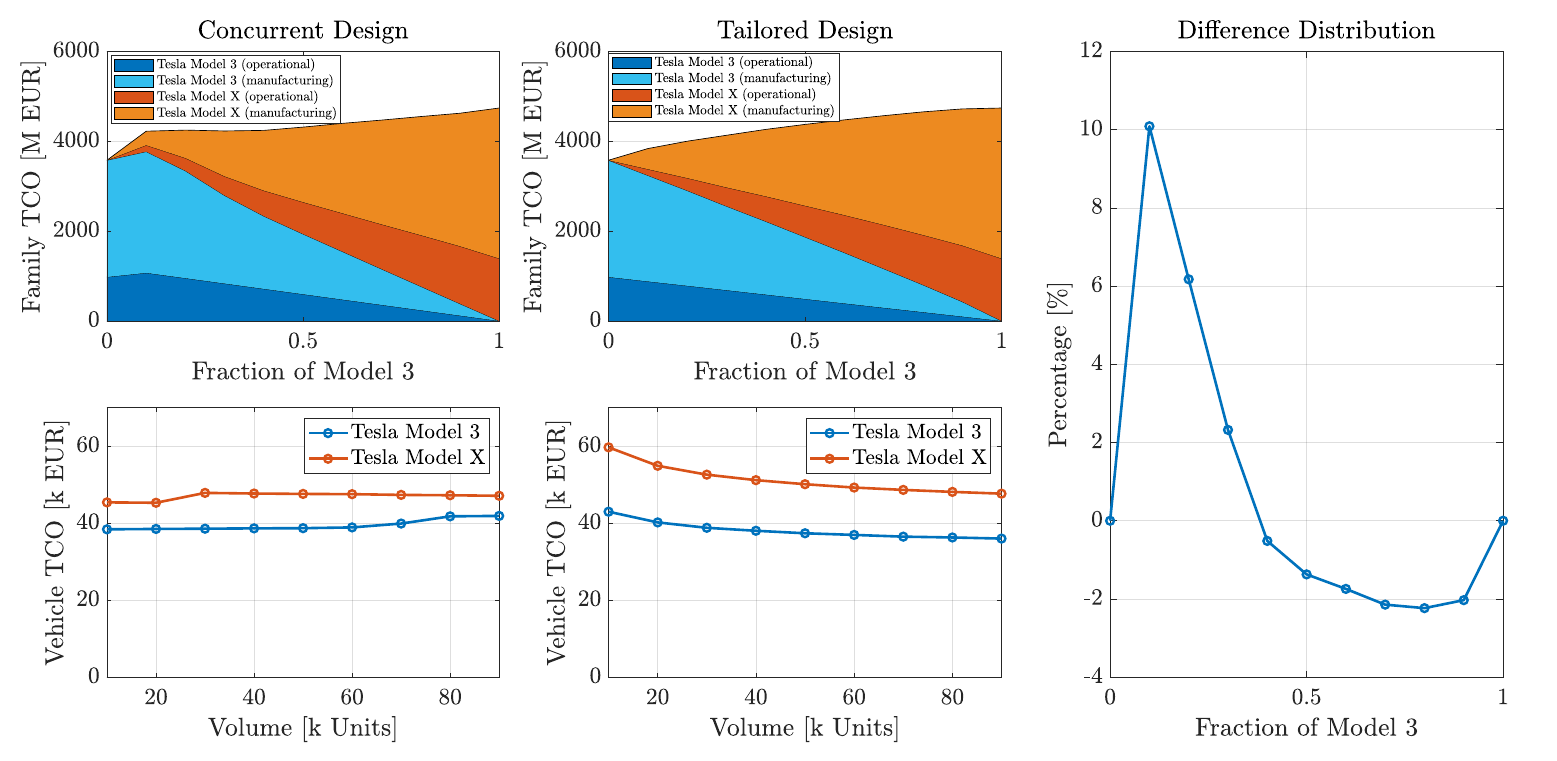}
	\caption{Sensitivity analysis of Tesla Model 3 and Tesla Model X with respect to the vehicle type fraction $w_i$ and a total production of 100000 units. For fractions 0 or 1, one of the vehicles is not produced while the other is optimized individually.}
	\label{fig:ratio2}
\end{figure*}	

\rwmargin{We analyze the sensitivity of the solution to the vehicle type fraction $w_i$, adjusting the ratio by 1/10 increments each time.
In the concurrent design optimization strategy, if the vehicles are similar enough to share the same powertrain, each vehicle's TCO depends on the total production volume, regardless of the specific production volume of each vehicle.
Conversely, using a vehicle-tailored design leads to different costs that depend on each vehicle's individual production volume.

As a matter of fact, in the scenario shown in Fig. C.16, the benefit of sharing a modular powertrain peaks around an approximately equal production fraction (though this depends on the respective costs of the vehicles), as an equal fraction does not significantly reduce the costs of any individual vehicle in a vehicle-tailored design.

When applying the concurrent design optimization strategy to vehicles that differ significantly and thus have different modularity of components (Fig. C.17), an uneven fraction could result in more expensive designs compared to a vehicle-tailored approach.
In fact, the smaller fraction could constrain the larger segment to adopt an oversized powertrain, increasing both component costs and energy consumption.}{R5:3}

\section*{Acknowledgments}
The authors would like to thank Ir. O.~J.~T.~Borsboom, F.~Vehlhaber, F.~Paparella, Dr.~D.~Herceg, and Dr.~I.~New for proofreading this paper.

\biboptions{sort&compress}
\bibliographystyle{unsrtnat}
\bibliography{../../../bibliography/main,../../../bibliography/SML_papers,references_eT}

\end{document}